\documentclass[11pt]{article} 
\begin{document} 
\newtheorem{prop}{Proposition}[section]
\newtheorem{Def}{Definition}[section] \newtheorem{theorem}{Theorem}[section]
\newtheorem{lemma}{Lemma}[section] \newtheorem{Cor}{Corollary}[section]

\title{\bf Well-posedness for a modified Zakharov system}
\author{{\bf Hartmut Pecher}\\
Fachbereich Mathematik und Naturwissenschaften\\
Bergische Universit\"at Wuppertal\\
Gau{\ss}str.  20\\
D-42097 Wuppertal\\
Germany\\
e-mail {\tt Hartmut.Pecher@math.uni-wuppertal.de}}
\date{}
\maketitle

\begin{abstract}
The Cauchy problem for a modified Zakharov system is proven to be locally 
well-posed for rough data in two and three space dimensions. In the three  
dimensional case the problem is globally well-posed for data with small energy. 
Under this assumption there also exists a global classical solution for 
sufficiently smooth data. 
\end{abstract}

\renewcommand{\thefootnote}{\fnsymbol{footnote}}
\footnotetext{\hspace{-1.8em}{\it 2000 Mathematics Subject Classification:} 
35Q55, 35L70 }
\normalsize 
\setcounter{section}{-1}
\section{Introduction}
The following system describes in plasma physics the nonlinear coupling of 
lower-hybrid waves, characterized by the complex amplitude $\varphi$ of the wave 
potential, with the much lower-frequency quasineutral density perturbations 
$\chi$ of the ion-acoustic type. It was introduced in \cite{S} as a variant of 
the standard Zakharov 
system which describes the phenomenon of Langmuir turbulence in a plasma. For 
details of the 
physical background and its derivation we refer to \cite{S}. The 
(2+1)-dimensional version reads as follows:
\begin{eqnarray}
\label{1.4}
i\frac{\partial}{\partial t} \Delta \varphi + \Delta^2 \varphi + \frac{1}{i} 
\nabla \varphi \cdot \overline{\nabla} \chi & = & 0 \\
\label{1.5}
\frac{\partial^2}{\partial t^2} \chi - \Delta \chi - \frac{1}{i} \Delta(\nabla 
\bar{\varphi} \cdot \overline{\nabla} \varphi) & = & 0 \, .
\end{eqnarray}
Here $\nabla$ denotes the usual gradient and $\overline{\nabla} = 
(\frac{\partial}{\partial x_2},-\frac{\partial}{\partial x_1})$ ,
and $\varphi$ and $\chi$ are respectively a complex-valued and a real-valued 
function defined for $(x,t) \in {\bf R}^2 \times {\bf R^+}$.

The initial conditions are
\begin{equation}
\label{1.3}
\varphi(x,0)= \varphi_0(x) \, , \, \chi(x,0)=\chi_0(x) \, , \, \frac{\partial 
\chi}{\partial t}(x,0)=\chi_1(x) \, .
\end{equation}
The functions $\varphi_0$ , $ \chi_0 $ , $ \chi_1$ are given in suitable 
Sobolev spaces.

A similar (3+1)-dimensional version of the Cauchy problem will also be 
considered, which reads as follows: 
\begin{eqnarray}
\label{1.1}
i\frac{\partial}{\partial t} \Delta \varphi + \Delta^2 \varphi + \frac{1}{i} 
(\nabla \varphi \times \nabla \chi) \cdot e & = & 0 \\
\label{1.2}
\frac{\partial^2}{\partial t^2} \chi - \Delta \chi - \frac{1}{i} \Delta(\nabla 
\bar{\varphi} \times \nabla \varphi) \cdot e & = & 0 \, .
\end{eqnarray}
Here $e$ is a constant vector in ${\bf R}^3$ and $\times$ denotes the vector 
product.

The most important question concerning the Cauchy problem is whether global 
smooth solutions exist for a class of smooth data. One way to attack this 
problem is to give a local well-posedness result for data with low regularity 
and then to use the conservation laws, especially the energy conservation, to 
extend this solution globally. It then remains to show that regular data lead to 
regular solutions. This program can in fact successfully be carried out, at 
least in 3+1 dimensions.

We are going to use the Fourier restriction norm method introduced by Bourgain 
\cite{B1},\cite{B2} to prove local existence and uniqueness of the problems 
also 
for rough data. It turns out that in 3+1 dimensions such a result is true
for the problem (\ref{1.1}),(\ref{1.2}),(\ref{1.3}) provided $B\varphi_0 \in 
H^k({\bf R}^3)$ , $ B\chi_0 \in H^l({\bf R}^3)$ , $B\chi_1 \in H^{l-1}({\bf 
R}^3)$ , where $B:= (-\Delta)^{\frac{1}{2}}$ , $ l \ge -1$ , $l+1 \le k \le 
l+2$ 
and $ k \ge \frac{l+2}{2} $ . So the lowest admissible pair is $(k,l) = 
(\frac{1}{2},-1)$ (cf. Theorem \ref{Theorem1}). It is also possible to treat the 
case $B\varphi_0 \in 
H^1({\bf R}^3) $ , $ \chi_0 \in L^2({\bf R}^3)$ , $B^{-1} \chi_1 \in L^2({\bf 
R}^3)$ . This is of particular interest, because in this case the conservation 
laws belonging to our problem (cf. (\ref{I_1}),(\ref{I_2}) below) can be used 
to 
give an a-priori bound for $\|B\varphi\|_{H^1} + \|\chi\|_{L^2} + 
\|B^{-1}\chi_t\|_{L^2}$, provided $\|B\varphi_0\|_{H^1} +\|\chi_0\|_{L^2} + 
\|B^{-1}\chi_1\|_{L^2}$ is sufficiently small. 
This allows to extend the solution globally in time, thus showing global 
well-posedness of the problem in energy space (Theorem \ref{Theorem2}).

It is also possible to refine these results in such a way (cf. Theorem 
\ref{Theorem3}) that one can show 
global well-posedness of the Cauchy problem for smoother data, especially 
proving the existence of global classical solutions under the above mentioned 
(weak) smallness assumption on the data (Theorem \ref{Theorem4}).

In 2+1 dimensions local well-posedness is proven for $B^{1+\epsilon} \varphi_0 
\in H^{k-\epsilon}({\bf R}^2)$ , $ B^{1-\delta} \chi_0 \in H^{l+\delta}({\bf 
R}^2)$ , $ B^{-\delta} \chi_1 \in H^{l+\delta}({\bf R}^2)$ , if $l \ge -1$ , 
$l+1 \le k \le l+2$ , $ k \ge \frac{l+2}{2} $ for $ 0 < \epsilon , \delta < 1 
$ (Theorem \ref{Theorem5})
. It is also possible to treat the case $B^{1+\epsilon} \varphi_0 \in 
H^{1-\epsilon}({\bf R}^2)$ , $\chi_0 \in L^2({\bf R}^2)$ , $ B^{-1}\chi_1 \in 
L^2({\bf R}^2) $ for $0 < \epsilon < 1 $ , but for global well-posedness one 
would need $\epsilon = 0$ , which is excluded here. The latter has to do with 
low frequency problems and the lack of a Sobolev embedding $\dot{H}^1 \subset 
L^{\infty}$ in two space dimensions.

This paper leaves open the question whether the results are optimal. In order to 
show the sharpness of the bilinear estimates one would need a number of 
counterexamples showing the necessity of the various conditions on the 
parameters involved. But even if this could be done this would not directly 
imply ill-posedness. A remarkable progress has been made in a recent paper by 
Holmer (\cite{H}) for the original Zakharov system in dimension 1+1, who made 
precise in which sense ill-posedness holds, if certain conditions on the 
parameters are violated. An idea could be to adapt these methods to the present 
more complicated higher dimensional situation, but I am not going to make such 
an attempt in this paper.

The technique of the proof relies on the pioneering works of Bourgain \cite{B1} 
and Kenig, Ponce and Vega \cite{KPV}, and especially on the paper of Ginibre - 
Tsutsumi - Velo 
\cite{GTV} for the corresponding problem for the original Zakharov system, which 
reads as follows:
\begin{eqnarray*}
i\frac{\partial}{\partial t} u + \Delta u  & = & nu \\
\frac{\partial^2}{\partial t^2} n - \Delta n & = & \Delta(|u|^2) \\
u(0) = u_0 \quad , \quad n(0) & = & n_0 \quad , \quad \frac{\partial n}{\partial 
t}(0) = n_1 \, .
\end{eqnarray*}
In 2+1 and 3+1 dimensions they showed local well-posedness for data $u_0 \in 
H^{k'},$  $n_0 \in H^{l'}$ , $n_1 \in H^{l'-1}$ under the assumptions $l' \ge 0 
$ , $ l' \le k' \le l'+1 $ , $ k' \ge \frac{l'+2}{2}.$ These conditions are in 
principle the same as ours (with $l' = l+1$ and $k'=k$), if one remarks that  
somehow $u $ can be identified with $(-\Delta)^{\frac{1}{2}} \varphi$ and $n$ 
with $\chi$ . Namely, after this identification and applying 
$(-\Delta)^{\frac{1}{2}}$ to the first equation of the Zakharov system we arrive 
at 
\begin{eqnarray*}
-i\frac{\partial}{\partial t} \Delta \varphi - \Delta^2 \varphi & = & 
(-\Delta)^{\frac{1}{2}}(\chi (-\Delta)^{\frac{1}{2}}\varphi) \\
\frac{\partial^2}{\partial t^2} \chi - \Delta \chi  & = & \Delta( 
|(-\Delta)^{\frac{1}{2}} \varphi|^2) \, ,
\end{eqnarray*}
which has a similar form as (\ref{1.1}),(\ref{1.2}) (just counting the number of 
derivatives), although the nonlinearities are of a different type.
 
Global well-posedness for the Zakharov system also holds for small data in 
two and three space dimensions \cite{BC}. A problem which is somehow related to 
the problem considered in the paper at hand has been treated in \cite{GY}. They 
however consider the 2-dimensional version with a weaker nonlinearity in the 
wave equation and prove global well-posedness for smooth data.

We will often use the notation $a+ = a+ \epsilon$ for a small $\epsilon > 0$ . 
Similarly, $a- = a- \epsilon$ and $a++ = a+2\epsilon$ .

The solution spaces are defined as follows: 
For $ k,l,b \in {\bf R} $ we denote by $X^{k,b}$ and $X_{\pm}^{l,b}$ the space 
such that $f \in {\cal S}'({\bf R}^n \times {\bf R})$ and 
$$ \|f\|_{X^{k,b}}^2 := \int \langle \tau + |\xi|^2 \rangle^{2b} \langle \xi 
\rangle^{2k} |\widehat{f}(\xi,\tau)|^2 \, d\xi d\tau < \infty $$
and
$$ \|f\|_{X_{\pm}^{l,b}}^2 := \int \langle \tau \pm |\xi| \rangle^{2b} \langle 
\xi \rangle^{2l} |\widehat{f}(\xi,\tau)|^2 \, d\xi d\tau < \infty \, , $$
respectively.
$ \dot{X}^{k,b}$ and $\dot{X}_{\pm}^{l,b}$ are defined by replacing $\langle 
\xi \rangle 
:= (1+|\xi|^2)^{\frac{1}{2}}$ by $|\xi|$ . $Y^k$ is defined with respect to 
$$ \|f\|_{Y^k} := \| \langle \tau + |\xi|^2 \rangle^{-1} \langle \xi \rangle^k 
\widehat{f}(\xi,\tau)\|_{L^2_{\xi}(L^1_{\tau})} $$
and $Y_{\pm}^l$ similarly by replacing $\langle \tau + |\xi|^2 \rangle^{-1}$ by 
$\langle \tau \pm |\xi| \rangle^{-1}$ . $\dot{Y}^k$ and $\dot{Y}_{\pm}^l$ are 
defined by replacing $\langle \xi \rangle$ by $|\xi|$ . We also use the 
corresponding restriction norm spaces $X^{k,b}[0,T]$ by its norm 
$\|f\|_{X^{k,b}[0,T]} := \inf_{\tilde{f}_{|[0,T]} = f} \| 
\tilde{f}\|_{X^{k,b}} 
$ and similarly the other cases.

We use the following standard facts about these spaces. Let $\psi$ denote a 
cut-off function in $C^{\infty}_0({\bf R})$ with $ supp \, \psi \subset 
(-2,2)$ 
, $ \psi = 1 $ on $ [-1,1]$ , $ \psi(t)=\psi(-t) $ , $ \psi(t) \ge 0$ , 
$\psi_{\delta}(t) := \psi(\frac{t}{\delta})$ , $ 0 < \delta \le 1 $ . Then the 
following estimates hold:
$$ \|\psi_{\delta} e^{it\Delta} f\|_{X^{k,b}} \le c \delta^{\frac{1}{2}-b} 
\|f\|_{H^k_x} \, , \, b \ge 0 $$
and similarly
$$ \|\psi_{\delta} e^{\pm itB} f\|_{X_{\pm}^{l,b}} \le c 
\delta^{\frac{1}{2}-b} 
\|f\|_{H^l_x} \, , \, b \ge 0 \, .$$
Moreover
\begin{equation}
\label{a}
\| \psi_{\delta} \int_0^t e^{-i(t-s)\Delta} f(s) \, ds \|_{X^{k,b}} \le c 
\delta^{1-b+b'} \|f\|_{X^{k,b'}} 
\end{equation}
for $b' \le 0 \le b \le b'+1$, $ b' > - \frac{1}{2} $ , $ \delta \le 1$ , and
\begin{equation}
\label{b}
\| \psi_{\delta} \int_0^t e^{-i(t-s)\Delta} f(s) \, ds \|_{X^{k,\frac{1}{2}}} 
\le c ( \|f\|_{X^{k,-\frac{1}{2}}} + \|f\|_{Y^k})
\end{equation}
as well as
\begin{equation}
\label{c}
\|\psi_{\delta} f \|_{X^{k,b}} \le c \delta^{-\epsilon} \|f\|_{X^{k,b}}
\end{equation}
for $ b \ge 0 $ , $ \epsilon > 0 $ .\\
Similar estimates hold for $X_{\pm}^{k,b}$, where $-\Delta$ is replaced by 
$B:=(-\Delta)^{\frac{1}{2}}$ .\\
Proofs can be found in \cite{GTV}.

The Strichartz estimates for the Schr\"odinger equation in ${\bf R}^n$ are 
given 
by
$$ \|e^{it\Delta} u_0\|_{L_t^q(L_x^r)} \le c \|u_0\|_{L_x^2} \, ,$$
if $ 0 \le \frac{2}{q} = n(\frac{1}{2}-\frac{1}{r}) < 1 $ . A direct 
consequence 
is (cf. \cite{GTV}, Lemma 2.4):
\begin{equation}
\label{0}
\|f\|_{L_t^q(L_x^r)} \le c \|f\|_{X^{0,b}} \, ,
\end{equation}
if $b_0 > \frac{1}{2}$, $ 0 \le b \le b_0 $ , $ 0\le \eta \le 1 $ , $ 
\frac{2}{q} = 1-\eta \frac{b}{b_0}$ , $ n(\frac{1}{2} - \frac{1}{r}) = 
(1-\eta)\frac{b}{b_0} $ . 

For the wave equation we only use
$$ \|e^{\pm itB} u_0\|_{L_t^{\infty}(L_x^2)} \le c \|u_0\|_{L_x^2} $$
and its consequence
\begin{equation}
\label{0'}
\|f\|_{L_t^q(L_x^2)} \le c \|f\|_{X_{\pm}^{0,b}} \, ,
\end{equation}
if $b_0 > \frac{1}{2} $ , $\frac{2}{q} = 1-\frac{b}{b_0} $ .

An important consequence for functions with a suitable support property is 
given 
by \cite{GTV}, Lemma 3.1, which we state as follows (for the Schr\"odinger 
equation):  
\begin{lemma}
\label{Lemma0}
Let $\sigma = \tau + |\xi|^2$ , $b_0 > \frac{1}{2}$ , $a \ge 0$ , $0\le \gamma 
\le 
1$ , $(1-\gamma)a \le b_0$ , $ a' \ge \gamma a.$ Define $\frac{2}{q} = 
1-\eta(1-\gamma)\frac{a}{b_0}$ , $ 
n(\frac{1}{2}-\frac{1}{r}):=(1-\eta)(1-\gamma)\frac{a}{b_0}$ . Let $v \in L^2$ 
be given such that ${\cal F}^{-1}(\langle \sigma \rangle^{-a'} \widehat{v})$ 
has 
support in $\{|t| \le cT\}$ . Then the following estimate holds:
$$ \|{\cal F}^{-1}(\langle \sigma \rangle^{-a} |\widehat{v}|)\|_{L_t^q(L_x^r)} 
\le c T^{\Theta} \|v\|_{L_x^2} \, ,$$
where $ \Theta = \gamma a(1-\frac{[a'-\frac{1}{2}]_+}{a'}) $ ,  
$[a'-\frac{1}{2}]_+ := a' - \frac{1}{2}$ , if $a' > \frac{1}{2}$ , $ := 
\epsilon$ , if $a'=\frac{1}{2}$ , $:=0,$ if $a' < \frac{1}{2}$ .
\end{lemma}
The proof is a combination of (\ref{0}), the support property and H\"older's 
inequality.\\
{\bf Remark:} 1. The same estimate is true for the wave equation with $\sigma 
:= 
\tau \pm |\xi|$ in the special case $\eta =1$ , $ r = 2 $ (by use of 
(\ref{0'})).\\
2. The statement of the Lemma without the factor $T^{\Theta}$ remains true, if 
no support property is assumed (with even a simpler proof).\\
For details we refer to \cite{GTV}.

{\bf Acknowledgment:} I am grateful to the referees for careful reading of the 
manuscript and helpful criticism.

\section{Conservation laws}
We now show that the system (\ref{1.1}),(\ref{1.2}) has two conserved 
quantities, namely
\begin{eqnarray}
\label{I_1}
I_1 & :=  & \int_{{\bf R}^3} |\nabla \varphi|^2 \, dx \\
\nonumber
I_2 & :=  & \int_{{\bf R}^3} |\Delta \varphi|^2 \, dx + \frac{1}{2} \int_{{\bf 
R}^3}(|(-\Delta)^{-\frac{1}{2}} \chi_t|^2 + |\chi|^2)dx + \frac{1}{i} 
\int_{{\bf 
R}^3}\hspace{-0.2em} \chi(\nabla \bar{\varphi} \times \nabla \varphi)\cdot e \, 
dx\\
\label{I_2}
& &
\end{eqnarray}
In order to show that $I_1$ is conserved we take the imaginary part of the 
scalar product of (\ref{1.1}) with $\varphi$. We use
\begin{eqnarray*}
\lefteqn{\Im \frac{1}{i} \langle(\nabla \varphi \times \nabla \chi)\cdot 
e,\varphi\rangle} \\
 & = & -\frac{1}{2} \int[(\varphi_{x_1} \chi_{x_2} - \varphi_{x_2} 
\chi_{x_1})\bar{\varphi} + \varphi(\bar{\varphi}_{x_1} \chi_{x_2} - 
\bar{\varphi}_{x_2} \chi_{x_1})]e_3 \, dx \\
 & & \mbox{+ 2 similar terms by permutation of the indices}
\end{eqnarray*} 
The first term is treated as follows
\begin{eqnarray*}
... & = & -\frac{e_3}{2} \int[(\varphi_{x_1} \chi)_{x_2} \bar{\varphi} - 
\varphi_{x_1 x_2} \chi \bar{\varphi} - (\varphi_{x_2}\chi)_{x_1} \bar{\varphi} 
+ 
\varphi_{x_2 x_1} \chi \bar{\varphi} \\
& & \qquad \; + \varphi(\bar{\varphi}_{x_1}\chi)_{x_2} - 
\varphi(\bar{\varphi}_{x_1 x_2} \chi) - \varphi(\bar{\varphi}_{x_2} 
\chi)_{x_1} 
+ \varphi \bar{\varphi}_{x_2 x_1} \chi] \, dx \\
& = & 0 \, .
\end{eqnarray*}
This implies that $I_1$ is conserved.

Next we show that $I_2$ is conserved. We take the real part of the scalar  
product of (\ref{1.1}) with $\varphi_t$ . We remark that
$$ \Re \langle i \Delta \varphi_t,\varphi_t\rangle = 0 \quad , \quad \Re 
\langle 
\Delta^2 \varphi,\varphi_t\rangle = \frac{1}{2} \frac{d}{dt} \|\Delta 
\varphi\|^2 $$
and
$$ \Re \frac{1}{i} \langle (\nabla \varphi \times \nabla \chi)\cdot 
e,\varphi_t\rangle = \frac{1}{2i}(\langle(\nabla \varphi \times \nabla 
\chi)\cdot e,\varphi_t\rangle - \langle(\nabla {\bar \varphi} \times \nabla 
\chi)\cdot e,\bar{\varphi}_t\rangle) \, .$$
Calculating $(\nabla \varphi \times \nabla \chi)\cdot e$ and taking its third 
term (the others are similar) we get
\begin{eqnarray*}
& & \hspace{-1.5em} \frac{e_3}{2i} \int((\varphi_{x_1} \chi_{x_2} - 
\varphi_{x_2} 
\chi_{x_1})\bar{\varphi}_t - (\bar{\varphi}_{x_1} \chi_{x_2} - 
\bar{\varphi}_{x_2} \chi_{x_1})\varphi_t) \, dx \\
 & & \hspace{-1em}= \frac{e_3}{2i} \int[(\varphi_{x_1} \chi)_{x_2} 
\bar{\varphi}_t - 
\varphi_{x_1 x_2} \chi \bar{\varphi}_t - (\varphi_{x_2}\chi)_{x_1} 
\bar{\varphi}_t + \varphi_{x_2 x_1} \chi \bar{\varphi}_t \\
 & & \hspace{-1em}\qquad -(\bar{\varphi}_{x_1} \chi)_{x_2} \varphi_t + 
\bar{\varphi}_{x_1 
x_2} \chi \varphi_t + (\bar{\varphi}_{x_2} \chi)_{x_1} \varphi_t - 
\bar{\varphi}_{x_2 x_1} \chi \varphi_t] \, dx \\
 & & \hspace{-1em}= \frac{e_3}{2i} \int (-\varphi_{x_1} \chi 
\bar{\varphi}_{tx_2} + 
\varphi_{x_2} \chi \bar{\varphi}_{tx_1} + \bar{\varphi}_{x_1} \chi 
\varphi_{tx_2} - \bar{\varphi}_{x_2} \chi \varphi_{tx_1}) \, dx \\
 & & \hspace{-1em}= \frac{e_3}{2i} \int \hspace{-0.2em}\chi(-\varphi_{x_1} 
\bar{\varphi}_{tx_2} + 
(\bar{\varphi}_{x_1} \varphi_{x_2})_t - \bar{\varphi}_{x_1} \varphi_{tx_2} + 
\bar{\varphi}_{x_1} \varphi_{tx_2} - (\bar{\varphi}_{x_2} \varphi_{x_1})_t + 
\bar{\varphi}_{tx_2} \varphi_{x_1}) dx \\
 & & \hspace{-1em}= \frac{e_3}{2i} \int \chi(\bar{\varphi}_{x_1} \varphi_{x_2} - 
\bar{\varphi}_{x_2} \varphi_{x_1})_t \, dx \, .
\end{eqnarray*} 
Thus we arrive at
\begin{eqnarray*}
\lefteqn{\Re \frac{1}{i} \langle (\nabla \varphi \times \nabla \chi)\cdot 
e,\varphi_t\rangle = \frac{1}{2i} \int \chi((\nabla \bar{\varphi} \times 
\nabla 
\varphi)\cdot e)_t \, dx} \\
& = & \frac{1}{2i} \frac{d}{dt} \int \chi(\nabla \bar{\varphi} \times \nabla 
\varphi)\cdot e \, dx - \frac{1}{2i} \int \chi_t (\nabla \bar{\varphi} \times 
\nabla{\varphi}) \cdot e \, dx \\
& = & \frac{1}{2i} \frac{d}{dt} \int \chi(\nabla \bar{\varphi} \times \nabla 
\varphi)\cdot e \, dx - \frac{1}{2} \int \chi_t (\Delta^{-1} \chi_{tt} - \chi) 
\, dx
\end{eqnarray*}
by using (\ref{1.2}). Now we have
\begin{eqnarray*}
- \frac{1}{2} \int \chi_t (\Delta^{-1} \chi_{tt} - \chi) \, dx & = & 
\frac{1}{2}(\langle (-\Delta)^{-\frac{1}{2}} \chi_t,(-\Delta)^{-\frac{1}{2}} 
\chi_{tt}\rangle + \langle \chi_t,\chi\rangle ) \\
& = & \frac{1}{4}\frac{d}{dt}(\|(-\Delta)^{-\frac{1}{2}} \chi_t\|^2 + 
\|\chi\|^2 
) \, .
\end{eqnarray*}
Summarizing we get
$$ \frac{d}{dt} \left( \|\Delta \varphi\|^2 + \frac{1}{2} 
(\|(-\Delta)^{-\frac{1}{2}} \chi_t\|^2+\|\chi\|^2) + \frac{1}{i} 
\int\chi(\nabla 
\bar{\varphi} \times \nabla \varphi) \cdot e \, dx \right) = 0 \, . $$

These two conservation laws imply an a-priori bound for the solution of our 
system (\ref{1.1}),(\ref{1.2}),(\ref{1.3}), provided suitable norms of the 
data 
are sufficiently small.
\begin{prop}
\label{Proposition}
Let $(\varphi,\chi)$ be a solution of (\ref{1.1}),(\ref{1.2}),(\ref{1.3}) with 
$B\varphi \in C^0([0,T],$ $ H^1({\bf R}^3))$ , $ \chi \in C^0([0,T],L^2({\bf 
R}^3))$ , $ B^{-1} \chi_t \in C^0([0,T],L^2({\bf R}^3))$ . Assume that the 
data 
fulfill
$$ \|B\varphi_0\|_{H^1} + \|\chi_0\|_{L^2} + \|B^{-1} \chi_t\|_{L^2} < 
\epsilon_0 $$
for a sufficiently small $\epsilon_0$ dependent only on the vector $e$ and 
some 
Sobolev embedding constants. Then for $t \in [0,T]$ :
$$ \|B\varphi(t)\|_{H^1} + \|\chi(t)\|_{L^2} + \|B^{-1} \chi_t(t)\|_{L^2} \le 
C_0 \, 
, 
$$
where $C_0$ is independent of $T$.
\end{prop}
{\bf Proof:} Consider the conserved quantity
$$ E(\varphi,\chi,\chi_t) := \|\Delta \varphi\|^2 + \frac{1}{2} \|\chi\|^2 + 
\frac{1}{2} \|B^{-1} \chi_t\|^2 + \frac{1}{i} \int \chi(\nabla \bar{\varphi} 
\times \nabla \varphi)\cdot e \, dx + \|\nabla \varphi\|^2 \, . $$
Now by the Sobolev embeddding $H^1({\bf R}^3) \subset L^4({\bf R}^3)$ :
\begin{eqnarray}
\nonumber
\frac{1}{2} \left|\int \chi(\nabla \bar{\varphi} \times \nabla \varphi) \times 
e 
\, dx\right| & \le & c \int |\chi| |\nabla \varphi|^2 \, dx \le \frac{1}{4} 
\int 
|\chi|^2 \, dx + c' \int |\nabla \varphi|^4 \, dx \\
\label{**}
& \le & \frac{1}{4} \int|\chi|^2 \, dx + c_0(\|\nabla \varphi\|^2 + \|\Delta 
\varphi\|^2)^2
\end{eqnarray} 
Defining
\begin{eqnarray*}
\lefteqn{\tilde{E}(\varphi_0,\chi_0,\chi_1)} \\ & := & \|\Delta \varphi_0\|^2 
+ 
\frac{1}{2} \|\chi_0\|^2 + \frac{1}{2} \|B^{-1} \chi_1\|^2 + \left| \int 
\chi_0(\nabla \bar{\varphi}_0 \times \nabla \varphi_0)\cdot e \, dx \right| + 
\|\nabla \varphi_0\|^2 \, , 
\end{eqnarray*}
we get
\begin{eqnarray*}
m(t) & := & \|\Delta \varphi\|^2 + \frac{1}{4} \|\chi\|^2 + \frac{1}{2} 
\|B^{-1} 
\chi_t \|^2 + \| \nabla \varphi \|^2 \\
& \le & \tilde{E}(\varphi_0,\chi_0,\chi_1) + c_0 (\|\nabla \varphi\|^2 + 
\|\Delta \varphi\|^2)^2 \, ,
\end{eqnarray*}
thus
$$ m(t) \le \tilde{E}(\varphi_0,\chi_0,\chi_1) + c_0 m(t)^2 \quad
\forall  t \in [0,T] \, .$$ 
Defining
$$ f(m) := \tilde{E}(\varphi_0,\chi_0,\chi_1) - m + c_0 m^2 $$
we get $ f(m(t)) \ge 0 $ $ \forall t \in [0,T]$ . $f$ has its only minimum in 
$m_0 = \frac{1}{2c_0}$ . For a suitably chosen $C_0$ our smallness assumption 
implies $ \tilde{E}(\varphi_0,\chi_0,\chi_1) < \frac{1}{4c_0} $ using 
(\ref{**}) 
above. This implies
$$ f(m_0) < \frac{1}{4c_0} - m_0 + c_0 m_0^2 = \frac{1}{4c_0} - \frac{1}{2c_0} 
+ 
c_0 \frac{1}{4c_0^2} = 0 \, . $$
Because $m(0) \le \tilde{E}(\varphi_0,\chi_0,\chi_1) < \frac{1}{4c_0} < m_0 $ 
and $f(m(0)) \ge 0$ ,  this implies $m(0) \le m_1$ , where $m_1$ is the smaller 
zero of $f(m)$ .
Because $m(t)$ is continuous and $f(m(t)) \ge 0$ we conclude $m(t) \le m_1$ 
$\forall t \in [0,T]$ and especially $m(t) \le m_0 $ $ \forall t \in [0,T]$ . 
Thus we have an a-priori bound for $m(t)$ , and the claim follows.

Concerning the (2+1)-dimensional problem the system 
(\ref{1.4}),(\ref{1.5}),(\ref{1.3}) has also two conserved quantities, namely
\begin{eqnarray*}
I_1 & := & \int_{{\bf R}^2} |\nabla \varphi|^2 \, dx \\
I_2 & := & \int_{{\bf R}^2} |\Delta \varphi|^2 \, dx + \frac{1}{2} \int_{{\bf 
R}^2}(|(-\Delta)^{-\frac{1}{2}} \chi_t|^2 + |\chi|^2)dx + \frac{1}{i} 
\int_{{\bf 
R}^2} \chi(\nabla \bar{\varphi} \cdot \overline{\nabla} \varphi) \, dx \, .\\
\end{eqnarray*}
This is shown in the same manner as in 3 dimensions. Moreover it is easy to see 
that these conservation laws imply an a-priori bound for $\|B\varphi\|_{H^1} + 
\|\chi\|_{L^2} + \|B^{-1} \chi_t\|_{L^2}$ , provided $\|B\varphi_0\|_{L^2}$ is 
sufficiently small. This follows immediately from a Gagliardo-Nirenberg type 
inequality for the cubic term in $I_2$, namely
\begin{eqnarray*}
\left| \int_{{\bf R}^2} \chi(\nabla \varphi \cdot \overline{\nabla} \varphi) 
\, 
dx \right| & \le & \frac{1}{4} \|\chi\|_{L^2}^2 + c \|\nabla \varphi\|_{L^2}^2 
\|\Delta \varphi\|_{L^2}^2 \\
& \le & \frac{1}{4} \|\chi\|_{L^2}^2 + \frac{1}{2} \|\Delta \varphi\|^2\, ,
\end{eqnarray*}
provided $ c\|\nabla \varphi_0\|_{L^2}^2 \le \frac{1}{2} $ .

The systems in 2+1 as well as in 3+1 dimensions can be transformed into a 
first 
order system in $t$ by defining
$$ \chi_{\pm} := \chi \pm i(-\Delta)^{-\frac{1}{2}}\frac{\partial 
\chi}{\partial 
t} \quad , \quad \chi = \frac{1}{2}(\chi_+ + \chi_-) \quad , \quad \chi_{\pm 
0} 
:= \chi_0 \pm i (-\Delta)^{-\frac{1}{2}} \chi_1 \, . $$
In 3+1 dimensions this leads to the system
\begin{eqnarray*}
i \frac{\partial}{\partial t} \Delta \varphi + \Delta^2 \varphi + 
\frac{1}{2i}(\nabla \varphi \times \nabla(\chi_+ + \chi_-)) \cdot e & = & 0 \\
i \frac{\partial}{\partial t} \chi_{\pm} \mp (\Delta)^{-\frac{1}{2}} 
\chi_{\pm} 
\pm \frac{1}{i} (\Delta)^{-\frac{1}{2}}(\nabla \bar{\varphi} \times \nabla 
\varphi) \cdot e & = & 0
\end{eqnarray*}
and
$$ \varphi(0) = \varphi_0 \quad , \quad \chi_{\pm}(0) = \chi_{\pm 0} \, .$$
The corresponding system of integral equations reads as follows:
\begin{eqnarray*}
(-\Delta)^{\frac{1}{2}} \varphi (t) & \hspace{-0.7em} = & 
\hspace{-0.7em}(-\Delta)^{\frac{1}{2}} e^{it\Delta} \varphi_0 - \frac{1}{2i} 
\int_0^t \hspace{-0.3em}e^{i(t-s)\Delta}(-\Delta)^{-\frac{1}{2}}((\nabla \varphi 
\times 
\nabla(\chi_+ + \chi_-)) \cdot e)  ds \\
(-\Delta)^{\frac{1}{2}} \chi_{\pm}(t) & \hspace{-0.7em}= & \hspace{-0.7em} 
(-\Delta)^{\frac{1}{2}}e^{\mp it(-\Delta)^{\frac{1}{2}}} \chi_{\pm 0} \mp 
\frac{1}{i} \int_0^t \hspace{-0.3em}  e^{\mp 
i(t-s)(-\Delta)^{\frac{1}{2}}}(-\Delta)((\nabla 
\bar{\varphi} \times \nabla \varphi) \cdot e)  ds  .
\end{eqnarray*}

\section{Local and global existence in 3+1 dimensions}
Concerning the system (\ref{1.1}),(\ref{1.2}),(\ref{1.3}), in order to prove 
local existence and uniqueness for solutions $B\varphi \in X^{k,b}[0,T]$ and 
$B\chi \in X^{l,b_1}_+[0,T] + X_-^{l,b_1}[0,T]$ we have to give estimates for 
the nonlinearities in spaces of the type $X^{k,b'}$ and $X_{\pm}^{l,b_1'}$ for 
some $b',b_1' \le 0$ , and in some limiting cases also in the spaces $Y^k$ and 
$Y_{\pm}^l$ , respectively, because in these cases we are forced to choose 
$b'= 
-\frac{1}{2}$ or $b_1' = -\frac{1}{2}$ (cf. (\ref{a}) and (\ref{b})).

In the sequel we use the notation 
$$ \xi:=\xi_1 - \xi_2 \, , \, \tau := \tau_1 - \tau_2 \, , \, \sigma_i := 
\tau_i 
+ |\xi_i|^2 \, (i=1,2) \, , \, \sigma:= \tau \pm |\xi| \, . $$
Then we have
\begin{equation}
\label{3.0}
|\xi_1|^2 - |\xi_2|^2 \mp |\xi| = \sigma_1 - \sigma_2 - \sigma \, . 
\end{equation}
Later we need the following elementary algebraic inequalities, which were 
essentially proven in (\cite{GTV}), Lemma 3.3. Here $\phi_E$ denotes the 
characteristic function of the set $E$.
\begin{lemma}
\begin{enumerate}
\item Let $y_1,y_2 \in {\bf R}$ and $z=y_1-y_2$ . Then for any $\lambda > 1$
\begin{equation}
\label{3.0'}
|z| \le \lambda |y_2| + \frac{\lambda}{\lambda -1} |y_1| 
\phi_{\{\frac{\lambda}{\lambda+1} \le \frac{|z|}{|y_1|} \le 
\frac{\lambda}{\lambda -1}\}} \, .
\end{equation}
\item Let $|\xi_1| \ge 2 |\xi_2|$. Then
\begin{eqnarray}
\label{3.1}
\langle \xi_1 \rangle^2 & \le & c (\langle \sigma \rangle + \langle \sigma_1 
\rangle + \langle \sigma_2 \rangle) \\
\label{3.2}
\langle \xi_1 \rangle^2 & \le & c(\langle \sigma \rangle + \langle \sigma_2 
\rangle + \langle \sigma_1 \rangle 
\phi_{\{c_1|\sigma_1| \le |\xi_1|^2 \le c_2|\sigma_1|\}}) \\
\label{3.3}
\langle \xi_1 \rangle^2 & \le & c(\langle \sigma_1\rangle + \langle \sigma_2 
\rangle + \langle \sigma \rangle 
\phi_{\{c_1|\sigma| \le |\xi|^2 \le c_2 |\sigma|\}}) \, ,
\end{eqnarray}
where $c,c_1,c_2 > 0$ .
\end{enumerate}
\end{lemma}
{\bf Proof:} (\ref{3.0'}) follows from the fact that $ \frac{\lambda 
-1}{\lambda} |z| \le |y_1| \le \frac{\lambda +1}{\lambda}|z| $ , if $|z| \ge 
\lambda |y_2|$ .\\
(\ref{3.1}) is implied by (\ref{3.0}) and the fact that $|\xi_1|^2 - |\xi_2|^2 
\mp |\xi| \sim |\xi_1|$ for large $|\xi_1|$ , and that $|\xi_1|^2 - |\xi_2|^2 
\mp |\xi|$ is bounded for small $|\xi_1|$ . \\
In order to prove (\ref{3.2}) we use (\ref{3.0'}) with $ z = |\xi_1|^2 - 
|\xi_2|^2 \mp |\xi| $ , and get for large $|\xi_1|$:
\begin{eqnarray*}
|\xi_1|^2 \sim ||\xi_1|^2 - |\xi_2|^2 \mp |\xi|| & \le & \lambda(|\sigma| + 
|\sigma_2|) + \frac{\lambda}{\lambda -1} |\sigma_1| 
\phi_{\{\frac{\lambda}{\lambda+1} \le \frac{||\xi_1|^2 - |\xi_2|^2 \mp 
|\xi||}{|\sigma_1|} \le \frac{\lambda}{\lambda -1} \}} \\
& \le & c(\langle \sigma \rangle + \langle \sigma_2 \rangle + \langle \sigma_1 
\rangle \phi_{\{c_1 |\sigma_1| \le |\xi_1|^2 \le c_2 |\sigma_2|\}}) \, .
\end{eqnarray*}
But (\ref{3.2}) is trivially also true for small $|\xi_1|$. \\
Finally, (\ref{3.3}) follows from (\ref{3.2}) by interchanging $\sigma$ and 
$\sigma_1$ and using $|\xi| \sim |\xi_1|$ .
\begin{lemma}
\label{LemmaA}
In space dimensions $n=2$ or $n=3$ let $ m > 0 $ , $ \frac{1}{2} \ge a,a_1,a_2 
\ge 0 $ satisfy $2(a+a_1+a_2)+m > \frac{n}{2} + 1$ and $a+a_1+a_2 > \frac{1}{2}$ 
. Let $v,v_1,v_2 \in 
L^2_{xt}$ 
be given such that ${\cal F}^{-1}(\langle \sigma \rangle^{-b} \widehat{v})$ 
and 
${\cal F}^{-1}(\langle \sigma_i \rangle^{-b_i} \widehat{v_i})$ are supported 
in 
$\{|t|\le cT\}$ for some $B \ge b \ge a$ , $B \ge b_i \ge a_i$ $(i=1,2)$. Then 
the following estimates hold with $\Theta = \Theta(a,a_1,a_2,m,B) > 0$ :
\begin{eqnarray*}
\int \frac{|\widehat{v}\widehat{v_1}\widehat{v_2}|}{\langle \sigma \rangle^a 
\langle \sigma_1 \rangle^{a_1} \langle \sigma_2 \rangle^{a_2} \langle \xi 
\rangle^m} & \le & c T^{\Theta} \|v\|_{L^2_{xt}} \|v_1\|_{L^2_{xt}} 
\|v_2\|_{L^2_{xt}} \, ,\\
\int \frac{|\widehat{v}\widehat{v_1}\widehat{v_2}|}{\langle \sigma \rangle^a 
\langle \sigma_1 \rangle^{a_1} \langle \sigma_2 \rangle^{a_2} \langle \xi_2 
\rangle^m} & \le & c T^{\Theta} \|v\|_{L^2_{xt}} \|v_1\|_{L^2_{xt}} 
\|v_2\|_{L^2_{xt}} \, .
\end{eqnarray*}
\end{lemma}
{\bf Remark:} Here and in the following integrals are always taken over 
$d\xi_1 \, d\xi_2 \, d\tau_1 \, d\tau_2$ and $\widehat{v} = 
\widehat{v}(\xi,\tau)$ , $\widehat{v_1} = \widehat{v_1}(\xi_1,\tau_1)$ 
,$\widehat{v_2} = \widehat{v_2}(\xi_2,\tau_2)$ . \\
{\bf Proof:} For the proof of the second inequality we refer to Lemma 
\ref{LemmaB} 
below. Just remark that we can assume $ m < \frac{n}{2} $ w.l.o.g. under our 
assumptions $2(a+a_1+a_2)+m > \frac{n}{2} +1$ and $ a+a_1+a_2 > \frac{1}{2}$ .

Next we prove the first inequality along the lines of \cite{GTV}, Lemma 3.2. 
We 
estimate using H\"older's inequality by
\begin{eqnarray}
\nonumber
c \|{\cal F}^{-1}(\langle \xi \rangle^{-m} \langle \sigma \rangle^{-a}  
|\widehat{v}|)\|_{L_t^q(L_x^r)} & \cdot & \|{\cal F}^{-1}( \langle \sigma_1 
\rangle^{-a_1} |\widehat{v_1}|)\|_{L_t^{q_1}(L_x^{r_1})} \\
\label{1} 
&  \cdot & \|{\cal F}^{-1}( \langle \sigma_2 \rangle^{-a_2} 
|\widehat{v_2}|)\|_{L_t^{q_2}(L_x^{r_2})}
\end{eqnarray}
with
\begin{eqnarray}
\label{2}
\frac{1}{q} + \frac{1}{q_1} + \frac{1}{q_2} & = & 1 \, ,\\ 
\label{3}
\frac{1}{r} + \frac{1}{r_1} + \frac{1}{r_2} & = & 1 \, .
\end{eqnarray}
Choose $b_0 = \frac{1}{2} +\epsilon$ , $\epsilon$ sufficiently small, and 
$0<\gamma,\eta<1$ such that
$$ \frac{2}{q_i} = 1 - \eta(1-\gamma) \frac{a_i}{b_0} \, (i=1,2) \, , \, 
\frac{2}{q} = 1 - (1-\gamma)\frac{a}{b_0} $$
(remark that $(1-\gamma)\max(a,a_1,a_2) < b_0$ , because $a,a_1,a_2 \le 
\frac{1}{2}$ , so that $q,q_1,q_2 \ge 2$). Now (\ref{2}) is equivalent to
\begin{equation}
\label{2'}
(1-\gamma)(a+\eta(a_1+a_2)) = b_0 \, .
\end{equation}
Concerning the x-integration we use the Sobolev embedding $H_x^{m,2} \subset 
L_x^r$ for
\begin{equation}
\label{4}
m > n(\frac{1}{2} - \frac{1}{r}) \ge 0
\end{equation}
and choose
\begin{equation}
\label{5}
n(\frac{1}{2}-\frac{1}{r_i}) = (1-\gamma)(1-\eta)\frac{a_i}{b_0} \, .
\end{equation} 
With these choices an application of Lemma \ref{Lemma0} (+ Remark 1) gives the 
desired bound. Now (\ref{3}) by use of (\ref{5}) reduces to
$$n(\frac{1}{2}-\frac{1}{r}) = n(\frac{1}{r_1} + \frac{1}{r_2} - \frac{1}{2}) 
= 
-n(\frac{1}{2} - \frac{1}{r_1}) - n(\frac{1}{2}-\frac{1}{r_2}) + \frac{n}{2} = 
\frac{n}{2} - (1-\gamma)(1-\eta)\frac{a_1+a_2}{b_0}  . $$
From (\ref{2'}) we get $ (1-\gamma)\frac{\eta(a_1+a_2)}{b_0} = 1 - (1-\gamma) 
\frac{a}{b_0} $ and thus $n(\frac{1}{2}-\frac{1}{r}) = 1 + \frac{n}{2} - 
(1-\gamma) \frac{a+a_1+a_2}{b_0}$ so that (\ref{4}) reduces to the condition
\begin{equation}
\label{6}
m > 1 + \frac{n}{2} - (1-\gamma)\frac{a+a_1+a_2}{b_0} \, .
\end{equation}
It remains to check (\ref{2'}) and (\ref{6}). (\ref{6}) can be fulfilled for a 
suitable $0<\gamma <1$ close to 0, if $b_0$ is close enough to $\frac{1}{2}$ 
under our assumption $2(a+a_1+a_2)+m > \frac{n}{2}+1$ . Concerning (\ref{2'}) 
we 
only remark that $(1-\gamma)a < \frac{1}{2} < b_0$ , whereas 
$(1-\gamma)(a+a_1+a_2) > b_0$ for small $\gamma >0$ and $b_0$ close to 
$\frac{1}{2}$ by the assumption $a+a_1+a_2 > \frac{1}{2}$ . So (\ref{2'}) can 
be 
fulfilled for a suitable $0 < \eta < 1$ . \\
{\bf Remark:} Lemma \ref{LemmaA} remains true, if one of the three factors 
does 
not fulfill the support property and at least one of the exponents $a,a_1,a_2$ 
belonging to the other two factors is strictly positive. This follows by using 
Remark 2 to Lemma \ref{Lemma0}.

We also need the following variant of the previous Lemma.
\begin{lemma}
\label{LemmaB}
In space dimensions $n=2$ or $n=3$ let $ \frac{n}{2} >m \ge 0 $ , $ 
\frac{1}{2} 
\ge a,a_1,a_2 \ge 0 $ , $a_1 > 0$ satisfy $2(a+a_1+a_2)+m > \frac{n}{2} + 1$. 
Let $v,v_1,v_2 \in L^2_{xt}$ be given such that ${\cal F}^{-1}(\langle \sigma 
\rangle^{-b} \widehat{v})$ and 
${\cal F}^{-1}(\langle \sigma_i \rangle^{-b_i} \widehat{v_i})$ are supported 
in 
$\{|t|\le cT\}$ for some $B \ge b \ge a$ , $B \ge b_i \ge a_i$ $(i=1,2)$. Then 
the following estimate holds with $\Theta = \Theta(a,a_1,a_2,m,B) > 0$ :
$$
\int \frac{|\widehat{v}\widehat{v_1}\widehat{v_2}|}{\langle \sigma \rangle^a 
\langle \sigma_1 \rangle^{a_1} \langle \sigma_2 \rangle^{a_2} | \xi_2|^m}  \le  
c T^{\Theta} \|v\|_{L^2_{xt}} \|v_1\|_{L^2_{xt}} \|v_2\|_{L^2_{xt}} \, . $$
\end{lemma} 
{\bf Proof:}
Again using a variant of the proof of \cite{GTV}, Lemma 3.2 we estimate the 
l.h.s. by H\"older's inequality as follows: 
\begin{eqnarray}
\nonumber
c \|{\cal F}^{-1}( \langle \sigma \rangle^{-a}  
|\widehat{v}|)\|_{L_t^q(L_x^2)} 
& \cdot & \|{\cal F}^{-1}( \langle \sigma_1 \rangle^{-a_1} 
|\widehat{v_1}|)\|_{L_t^{q_1}(L_x^{r_1})} \\
\label{1*} 
&  \cdot & \|{\cal F}^{-1}(|\xi_2|^{-m} \langle \sigma_2 \rangle^{-a_2} 
|\widehat{v_2}|)\|_{L_t^{q_2}(L_x^{r_2})}
\end{eqnarray}
with
\begin{eqnarray}
\label{2*}
\frac{1}{q} + \frac{1}{q_1} + \frac{1}{q_2} & = & 1 \, ,\\ 
\label{3*}
 \frac{1}{r_1} + \frac{1}{r_2} & = & \frac{1}{2} \, . 
\end{eqnarray}
Choose $b_0 = \frac{1}{2} +\epsilon$ , $\epsilon$ sufficiently small, and 
$0<\gamma,\eta<1$ such that
$$ \frac{2}{q_i} = 1 - \eta(1-\gamma) \frac{a_i}{b_0} \, (i=1,2) \, , \, 
\frac{2}{q} = 1 - (1-\gamma)\frac{a}{b_0} $$
(remark that $(1-\gamma)\max(a,a_1,a_2) < b_0$ , because $a,a_1,a_2 \le 
\frac{1}{2}$ , so that $q,q_1,q_2 \ge 2$). Now (\ref{2*}) is equivalent to
\begin{equation}
\label{2'*}
(1-\gamma)(a+\eta(a_1+a_2)) = b_0 \, .
\end{equation}
Concerning the x-integration we use the Sobolev embedding $\dot{H}_x^{m,r_2'} 
\subset L_x^{r_2}$ provided
\begin{equation}
\label{4*}
m = n(\frac{1}{r_2'} - \frac{1}{r_2}) \ge 0
\end{equation} 
and $r_2 \neq \infty$ . This last condition is by (\ref{3*}) equivalent to 
$r_1 
\neq 2$ . We now choose $r_1$ such that
\begin{equation}
\label{5*}
n(\frac{1}{2}-\frac{1}{r_1}) := (1-\gamma)(1-\eta) \frac{a_1}{b_0} \, .
\end{equation}
This is strictly positive, because $a_1 > 0$ . Thus $r_1 \neq 2$ and $r_2 \neq 
\infty$ is fulfilled. Now we choose $r_2'$ such that
\begin{equation}
\label{6*}
n(\frac{1}{2}-\frac{1}{r_2'}) := (1-\gamma)(1-\eta) \frac{a_2}{b_0} \, .
\end{equation}
With these choices we can estimate (\ref{1*}) by $ \,cT^{\Theta} 
\|v\|_{L^2_{xt}}
\|v_1\|_{L^2_{xt}} \|v_2\|_{L^2_{xt}} \, $ using Lemma \ref{Lemma0} (+ Remark 
1).
Now we compute using (\ref{3*}),(\ref{5*}),(\ref{6*}):
\begin{eqnarray*}
n(\frac{1}{r_2'}-\frac{1}{r_2}) & = & n(\frac{1}{r_2'} + \frac{1}{r_1} - 
\frac{1}{2}) = \frac{n}{2} - (1-\gamma)(1-\eta)\frac{a_1+a_2}{b_0} \\ 
& = & \frac{n}{2} - (1-\gamma)\frac{a_1+a_2}{b_0} + 
\eta(1-\gamma)\frac{a_1+a_2}{b_0} \, .
\end{eqnarray*}
From (\ref{2'*}) we get $ (1-\gamma)\frac{\eta(a_1+a_2)}{b_0} = 1 - 
(1-\gamma)\frac{a}{b_0} $ and thus
$$ n(\frac{1}{r_2'} - \frac{1}{r_2}) = 1 + \frac{n}{2} - (1-\gamma) 
\frac{a+a_1+a_2}{b_0} \, . $$
Thus (\ref{4*}) reduces to
\begin{equation}
\label{7*}
m = 1 + \frac{n}{2} - (1-\gamma)\frac{a+a_1+a_2}{b_0} \Longleftrightarrow 
(1-\gamma)(a+a_1+a_2) = b_0(\frac{n}{2}+1-m) \, .
\end{equation} 
It remains to fulfill (\ref{2'*}) and (\ref{7*}). (\ref{7*}) can be fulfilled 
with a suitable $0 < \gamma < 1$, if $b_0$ is close enough to $\frac{1}{2}$ 
under our assumption $2(a+a_1+a_2)+m > \frac{n}{2}+1$ . It remains to fulfill 
(\ref{2'*}). By (\ref{7*}) and $ m < \frac{n}{2} $ we have $ 
(1-\gamma)(a+a_1+a_2) > b_0 $ , whereas $(1-\gamma)a < \frac{1}{2} < b_0$ , so 
that (\ref{2'*}) can be fulfilled by a suitable choice of $\eta \in (0,1)$ . 
\\
{\bf Remark:} Similarly as for Lemma \ref{LemmaA} it is sufficient here to 
have 
the support property for only two of the three factors, provided at least one 
of 
the exponents $a,a_1,a_2$ belonging to the other two factors is strictly 
positive.

In the following $D$ denotes any first order spatial derivative.
\begin{lemma}
\label{LemmaC}
In space dimension $n=3$ assume $ l \ge -1 $ , $ k \ge l+1 $ , $ k < l+2 $ 
with the exception of $(k,l)=(0,-1)$ . $\varphi$ and $\chi$ are given with 
support in $\{|t| \le cT\}$ . Then the following estimate holds:
$$ \|(-\Delta)^{-\frac{1}{2}} (D\varphi D\chi)\|_{X^{k,-\frac{1}{2}+}} \le 
cT^{\Theta} \|D\varphi\|_{X^{k,\frac{1}{2}}} 
\|D\chi\|_{X_{\pm}^{l,\frac{1}{2}}} $$
with $ \Theta = \Theta(k,l) > 0 $ .
\end{lemma} 
{\bf Remark:} Trivially we can replace $\|D\chi\|_{X_{\pm}^{l,\frac{1}{2}}}$ 
by $\|\chi\|_{\dot{H}_{\pm}^{l+1,\frac{1}{2}}}$ , if $ l \le 0 $ .\\
{\bf Proof:} Defining $\widehat{v} := \langle \xi \rangle^l \langle \sigma 
\rangle^{\frac{1}{2}} \widehat{D\chi}$ , $ \widehat{v_2}:= \langle \xi_2  
\rangle^k \langle \sigma_2 \rangle^{\frac{1}{2}} \widehat{D\varphi} $ and 
$\widehat{\psi} := \langle \xi_1 \rangle^k \langle \sigma_1 
\rangle^{-\frac{1}{2}+} \widehat{v_1}, $ where $ v_1 \in L^2_{xt} $ , we have 
$ \|v\|_{L^2_{xt}} = \|D\chi\|_{X_{\pm}^{l,\frac{1}{2}}} $ , $ 
\|v_2\|_{L^2_{xt}} = \|D\varphi\|_{X^{k,\frac{1}{2}}} $ and $ 
\|v_1\|_{L^2_{xt}} = \|\psi\|_{X^{-k,\frac{1}{2}-}} $ . This generic function 
$\psi$ in $X^{-k,\frac{1}{2}-}$ can be assumed to have support in $\{|t| \le 
cT\}$ , too. Thus we have: the support of ${\cal F}^{-1}(\langle \sigma 
\rangle^{-\frac{1}{2}} \widehat{v})$ , ${\cal F}^{-1}(\langle \sigma_2 
\rangle^{-\frac{1}{2}} \widehat{v_2})$ and ${\cal F}^{-1}(\langle \sigma_1 
\rangle^{-\frac{1}{2}+} \widehat{v_1})$ is contained in $\{|t| \le cT\}$ . We 
thus have to show:
\begin{equation}
\label{*}
S:=\left| \int \frac{\widehat{v} \widehat{v_1} \widehat{v_2} |\xi_1|^{-1} 
\langle \xi_1 \rangle^k}{\langle \xi \rangle^l \langle \xi_2 \rangle^k \langle 
\sigma \rangle^{\frac{1}{2}} \langle \sigma_1 \rangle^{\frac{1}{2}-} \langle 
\sigma_2 \rangle^{\frac{1}{2}}} \right| \le cT^{\Theta} \|v\|_{L^2_{xt}} 
\|v_1\|_{L^2_{xt}} \|v_2\|_{L^2_{xt}} \, .
\end{equation}  
{\bf Region A:} $|\xi_1| \le \frac{1}{2} |\xi_2|$.\\
In this case we have $|\xi| \sim|\xi_2|$ , thus
$$ S \le c \int \frac{|\widehat{v} \widehat{v_1} \widehat{v_2}| |\xi_1|^{-1} 
\langle \xi_1 \rangle^k}{ \langle \xi_2 \rangle^{k+l} \langle \sigma 
\rangle^{\frac{1}{2}} \langle \sigma_1 \rangle^{\frac{1}{2}-} \langle \sigma_2 
\rangle^{\frac{1}{2}}} \, . $$ 
{\bf Case 1:} $ k < 1 $ , $ k+l \le 0 $ . \\
We use the estimate (cf. (\ref{3.1})) $\langle \xi_2 \rangle \le (\langle 
\sigma \rangle + \langle \sigma_1 \rangle + \langle \sigma_2 
\rangle)^{\frac{1}{2}} $ and get
$$ S \le c \int \frac{|\widehat{v} \widehat{v_1} \widehat{v_2}| (\langle 
\sigma \rangle + \langle \sigma_1 \rangle + \langle \sigma_2 
\rangle)^{\frac{-k-l}{2}}}{|\xi_1| \langle \xi_1 \rangle^{-k} \langle \sigma 
\rangle^{\frac{1}{2}} \langle \sigma_1 \rangle^{\frac{1}{2}-} \langle \sigma_2 
\rangle^{\frac{1}{2}}} \, . $$
Because under our assumptions $-k-l < 1$ , we get three terms with positive 
powers of the $\sigma$ - modules in the denominator. \\
{\bf a.} We consider first the case $|\xi_1| \ge 1$ , where we have
$$ S \le c \int \frac{|\widehat{v} \widehat{v_1} \widehat{v_2}| (\langle 
\sigma \rangle + \langle \sigma_1 \rangle + \langle \sigma_2 
\rangle)^{\frac{-k-l}{2}}}{ \langle \xi_1 \rangle^{1-k} \langle \sigma 
\rangle^{\frac{1}{2}} \langle \sigma_1 \rangle^{\frac{1}{2}-} \langle \sigma_2 
\rangle^{\frac{1}{2}}} \, . $$
We use Lemma \ref{LemmaA} with e.g. $a=\frac{1}{2}+\frac{k+l}{2}$, $ a_1 = 
\frac{1}{2}-$ , $a_2 =\frac{1}{2}$ , $ m=1-k$ (and similar choices in the 
other cases) and get $2(a+a_1+a_2)+m = l+4- > \frac{5}{2}$ for $ l > 
-\frac{3}{2}$ , $a+a_1+a_2 = \frac{3}{2} + \frac{k+l}{2} - > \frac{1}{2} $ and 
$a,a_1,a_2 \le \frac{1}{2}$, because $k+l \le 0$ .\\
{\bf b.} In the case $|\xi_1| \le 1$ we get
$$ S \le c \int \frac{|\widehat{v} \widehat{v_1} \widehat{v_2}| (\langle 
\sigma \rangle + \langle \sigma_1 \rangle + \langle \sigma_2 
\rangle)^{\frac{-k-l}{2}}}{|\xi_1| \langle \sigma \rangle^{\frac{1}{2}} 
\langle \sigma_1 \rangle^{\frac{1}{2}-} \langle \sigma_2 
\rangle^{\frac{1}{2}}} \, . $$
Similarly as before we use Lemma \ref{LemmaB} with $m=1$ and get 
$2(a+a_1+a_2)+m=k+l+4 \ge -1+4 =3$ , thus the desired estimate. \\
{\bf Case 2:} $ k \le 1$ , $ k+l > 0$ . \\
We get
$$ S \le c \int \frac{|\widehat{v} \widehat{v_1} \widehat{v_2}| |\xi_1|^{-1} 
\langle \xi_1 \rangle^k}{ \langle \xi_1 \rangle^{k+l} \langle \sigma 
\rangle^{\frac{1}{2}} \langle \sigma_1 \rangle^{\frac{1}{2}-} \langle \sigma_2 
\rangle^{\frac{1}{2}}} =  c  \int \frac{|\widehat{v} \widehat{v_1} 
\widehat{v_2}|}{|\xi_1| \langle \xi_1 \rangle^{l} \langle \sigma 
\rangle^{\frac{1}{2}} \langle \sigma_1 \rangle^{\frac{1}{2}-} \langle \sigma_2 
\rangle^{\frac{1}{2}}} \, . $$ 
{\bf a.} $|\xi_1| \ge 1 $ . \\
By $ l \ge -1 $ we get
$$ S \le c  \int \frac{|\widehat{v} \widehat{v_1} \widehat{v_2}|}{\langle 
\sigma \rangle^{\frac{1}{2}} \langle \sigma_1 \rangle^{\frac{1}{2}-} \langle 
\sigma_2 \rangle^{\frac{1}{2}}} \, . $$
This can be handled by Lemma \ref{LemmaB} with $a=a_2=\frac{1}{2} $ , $ 
a_1=\frac{1}{2}-$ , $m=0$ . \\
{\bf b.} $|\xi_1| \le 1$ . \\
$$ S \le c  \int \frac{|\widehat{v} \widehat{v_1} \widehat{v_2}|}{|\xi_1| 
\langle \sigma \rangle^{\frac{1}{2}} \langle \sigma_1 \rangle^{\frac{1}{2}-} 
\langle \sigma_2 \rangle^{\frac{1}{2}}} \, . $$
We use Lemma \ref{LemmaB} with $a=a_2=\frac{1}{2}$ , $ a_1=\frac{1}{2}-$ , 
$m=1$ .\\
{\bf Case 3:} $ k \ge 1 $ . \\
{\bf a.} $|\xi_1| \ge 1$ . \\
Using $|\xi_1| \le \frac{1}{2}|\xi_2|$ and $l \ge -1$ we get
\begin{eqnarray*}
S  \le c \int \frac{|\widehat{v} \widehat{v_1} \widehat{v_2}|\langle \xi_1 
\rangle^{k-1}}{ \langle \xi_2 \rangle^{k+l} \langle \sigma 
\rangle^{\frac{1}{2}} \langle \sigma_1 \rangle^{\frac{1}{2}-} \langle \sigma_2 
\rangle^{\frac{1}{2}}} & \le & c \int \frac{|\widehat{v} \widehat{v_1} 
\widehat{v_2}|}{ \langle \xi_2 \rangle^{l+1} \langle \sigma 
\rangle^{\frac{1}{2}} \langle \sigma_1 \rangle^{\frac{1}{2}-} \langle \sigma_2 
\rangle^{\frac{1}{2}}} \\
& \le & c \int \frac{|\widehat{v} \widehat{v_1} \widehat{v_2}|}{ \langle 
\sigma \rangle^{\frac{1}{2}} \langle \sigma_1 \rangle^{\frac{1}{2}-} \langle 
\sigma_2 \rangle^{\frac{1}{2}}} \, .
\end{eqnarray*}
This can be handled by Lemma \ref{LemmaB} with $a=a_2=\frac{1}{2} $ , $ a_1 = 
\frac{1}{2}- $ , $ m=0 $ . \\
{\bf b.} $ |\xi_1| \le 1 $ . \\
Using $k+l \ge 1+l \ge 0$ we get
$$S  \le c \int \frac{|\widehat{v} \widehat{v_1} \widehat{v_2}|}{|\xi_1| 
\langle \xi_2 \rangle^{k+l} \langle \sigma \rangle^{\frac{1}{2}} \langle 
\sigma_1 \rangle^{\frac{1}{2}-} \langle \sigma_2 \rangle^{\frac{1}{2}}} 
 \le c \int \frac{|\widehat{v} \widehat{v_1} \widehat{v_2}|}{|\xi_1| \langle 
\sigma \rangle^{\frac{1}{2}} \langle \sigma_1 \rangle^{\frac{1}{2}-} \langle 
\sigma_2 \rangle^{\frac{1}{2}}} \, . $$
Now we use Lemma \ref{LemmaB} with $a=a_2=\frac{1}{2}$ , $ a_1 = \frac{1}{2}-$ 
, $m=1$ .\\
{\bf Region B:} $\frac{1}{2} |\xi_2| \le |\xi_1| \le 2|\xi_2|$ ($\Rightarrow \, 
|\xi| \le 3|\xi_1|,3|\xi_2|$). \\
We have
$$ S \le c \int \frac{|\widehat{v} \widehat{v_1} \widehat{v_2}| \langle \xi 
\rangle^{-l}}{|\xi_1| \langle \sigma \rangle^{\frac{1}{2}} \langle \sigma_1 
\rangle^{\frac{1}{2}-} \langle \sigma_2 \rangle^{\frac{1}{2}}} \, . $$
If $l \ge 0$ we arrive at the same integral as in Region A, Case 3b. \\
If $ -1 \le l <0$ we estimate as follows:
$$ S \le c \int \frac{|\widehat{v} \widehat{v_1} \widehat{v_2}| \langle \xi_1 
\rangle^{-l}}{|\xi_1| \langle \sigma \rangle^{\frac{1}{2}} \langle \sigma_1 
\rangle^{\frac{1}{2}-} \langle \sigma_2 \rangle^{\frac{1}{2}}} \le c \int 
\frac{|\widehat{v} \widehat{v_1} \widehat{v_2}| \langle \xi_1 \rangle}{|\xi_1| 
\langle \sigma \rangle^{\frac{1}{2}} \langle \sigma_1 \rangle^{\frac{1}{2}-} 
\langle \sigma_2 \rangle^{\frac{1}{2}}} \, . $$
In the case $|\xi_1| \le 1$ and $|\xi_1| \ge 1$ we arrive at the same integral 
as in Region A, Case 3b and Case 3a, respectively.\\
{\bf Region C:} $|\xi_1| \ge 2 |\xi_2| $ ($\Rightarrow \, |\xi| \sim |\xi_1|$). 
\\
We get 
$$ S \le c \int \frac{|\widehat{v} \widehat{v_1} \widehat{v_2}| |\xi_1|^{-1} 
\langle \xi_1 \rangle^{k-l}}{\langle \xi_2 \rangle^k  \langle \sigma 
\rangle^{\frac{1}{2}} \langle \sigma_1 \rangle^{\frac{1}{2}-} \langle \sigma_2 
\rangle^{\frac{1}{2}}} \, . $$
{\bf a.} $|\xi_1| \le 1 $ . \\
This implies $|\xi_2| \le \frac{1}{2}$ , so that we again arrive at the same 
term as in Region A, Case 3b. \\
{\bf b.} $|\xi_1| \ge 1$ . \\
Because $k \ge l+1$ by assumption, we get by (\ref{3.1}) : 
$$ S \le c \int \frac{|\widehat{v} \widehat{v_1} \widehat{v_2}| \langle \xi_1 
\rangle^{k-l-1}}{\langle \xi_2 \rangle^k  \langle \sigma \rangle^{\frac{1}{2}} 
\langle \sigma_1 \rangle^{\frac{1}{2}-} \langle \sigma_2 \rangle^{\frac{1}{2}}} 
\le c \int \frac{|\widehat{v} \widehat{v_1} \widehat{v_2}| (\langle \sigma 
\rangle + \langle \sigma_1 \rangle + \langle \sigma_2 
\rangle)^{\frac{k-l-1}{2}}}{\langle \xi_2 \rangle^k  \langle \sigma 
\rangle^{\frac{1}{2}} \langle \sigma_1 \rangle^{\frac{1}{2}-} \langle \sigma_2 
\rangle^{\frac{1}{2}}} \, . $$
We remark that our assumption $k < l+2$ implies that the exponents of the 
$\sigma$ - modules in the denominator are positive. Using Lemma \ref{LemmaA} 
with e.g. $a=\frac{1}{2}=a_2,$  $ a_1 = \frac{1}{2} - \frac{k-l-1}{2}-$ , $ m=k 
> 
0 $ , thus $2(a+a_1+a_2)+m=4+l- > \frac{5}{2}$ for $l>-\frac{3}{2}$, we get the 
desired bound. 
\begin{Cor}
\label{Corollary to Lemma C}
Under the assumptions of Lemma \ref{LemmaC} we have for $k \ge 1$:
$$\|(-\Delta)^{-\frac{1}{2}} (D\varphi D\chi)\|_{X^{k,-\frac{1}{2}+}} \le 
cT^{\Theta}(\|D\varphi\|_{X^{1,\frac{1}{2}}} 
\|\chi\|_{X_{\pm}^{l+1,\frac{1}{2}}} + \|D\varphi\|_{X^{k,\frac{1}{2}}} 
\|\chi\|_{X_{\pm}^{0,\frac{1}{2}}}) \, . $$
\end{Cor} 
{\bf Proof:} We use Lemma \ref{LemmaC} with $k=1-$ , $ l=-1 $ :
\begin{eqnarray*}
\|(-\Delta)^{-\frac{1}{2}} (D\varphi D\chi)\|_{X^{1-,-\frac{1}{2}+}} & \le & 
cT^{\Theta}\|D\varphi\|_{X^{1-,\frac{1}{2}}} 
\|D\chi\|_{X_{\pm}^{-1,\frac{1}{2}}} \\
& \le & cT^{\Theta} \|D\varphi\|_{X^{1-,\frac{1}{2}}} 
\|\chi\|_{X_{\pm}^{0,\frac{1}{2}}} \, . 
\end{eqnarray*}
Applying the elementary inequality $\langle \xi_1 \rangle^{k-1+} \le c(\langle 
\xi \rangle^{k-1+} + \langle \xi_2 \rangle^{k-1+})$ in the Fourier variables we 
arrive at
\begin{eqnarray*}
\|(-\Delta)^{-\frac{1}{2}} (D\varphi D\chi)\|_{X^{k,-\frac{1}{2}+}} & 
\hspace{-0.3em} \le & \hspace{-0.3em} 
cT^{\Theta}(\|D\varphi\|_{X^{1-,\frac{1}{2}}} 
\|\chi\|_{X_{\pm}^{k-1+,\frac{1}{2}}} + \|D\varphi\|_{X^{k,\frac{1}{2}}} 
\|\chi\|_{X_{\pm}^{0,\frac{1}{2}}}) \\
& \hspace{-0.3em} \le & \hspace{-0.3em} 
cT^{\Theta}(\|D\varphi\|_{X^{1,\frac{1}{2}}} 
\|\chi\|_{X_{\pm}^{l+1,\frac{1}{2}}} + \|D\varphi\|_{X^{k,\frac{1}{2}}} 
\|\chi\|_{X_{\pm}^{0,\frac{1}{2}}}) \, . 
\end{eqnarray*} 
\begin{lemma}
\label{LemmaD}
In space dimension $n=3$ assume $l \ge -1$ , $ k \ge \frac{l+2}{2} $ , $ k > l+1 
$ , and let $\varphi_1$ , $ \varphi_2 $ be supported in $\{|t| \le cT\}$ . Then 
the following estimate holds:
$$ \|D\bar{\varphi}_1 D \varphi_2\|_{X_{\pm}^{l+2,-\frac{1}{2}+}} \le c 
T^{\Theta} \|D \varphi_1\|_{X^{k,\frac{1}{2}}} \|D 
\varphi_2\|_{X^{k,\frac{1}{2}}} $$
with $\Theta = \Theta(k,l) > 0$ .
\end{lemma}
{\bf Remark:} Trivially we can replace $X_{\pm}^{l+2,-\frac{1}{2}+}$ by 
$\dot{X}_{\pm}^{l+2,-\frac{1}{2}+}$ . \\
{\bf Proof:} Defining $\,  \widehat{v_1} := \langle \xi_1 \rangle^k \langle 
\sigma_1 \rangle^{\frac{1}{2}} \widehat{D\varphi_1} \, $ , $ \, \widehat{v_2} := 
\langle \xi_2 \rangle^k \langle \sigma_2 \rangle^{\frac{1}{2}} 
\widehat{D\varphi_2}\,  $ and $\, \widehat{\psi} := \langle \xi \rangle^{l+2}$ $ 
\langle \sigma \rangle^{-\frac{1}{2}+} \widehat{v} \,$ , where $ v \in L^2 $ , 
we have to show 
$$  W = \left| \int \frac{\widehat{v} \widehat{v_1} \widehat{v_2} \langle \xi 
\rangle^{l+2}}{\langle \xi_1 \rangle^k \langle \xi_2 \rangle^k  \langle \sigma 
\rangle^{\frac{1}{2}-} \langle \sigma_1 \rangle^{\frac{1}{2}} \langle \sigma_2 
\rangle^{\frac{1}{2}}}\right| \le c T^{\Theta} \|v\|_{L^2_{xt}} 
\|v_1\|_{L^2_{xt}} \|v_2\|_{L^2_{xt}} \, . $$ 
{\bf Region A:} $\frac{|\xi_2|}{2} \le |\xi_1| \le 2|\xi_2|$ $( \Rightarrow \, 
|\xi| \le 3|\xi_1|,3|\xi_2|)$ . \\
This gives
$$  W \le c \int \frac{|\widehat{v} \widehat{v_1} \widehat{v_2}| \langle \xi_1 
\rangle^{l+2-2k}}{\langle \sigma \rangle^{\frac{1}{2}-} \langle \sigma_1 
\rangle^{\frac{1}{2}} \langle \sigma_2 \rangle^{\frac{1}{2}}} \le c \int 
\frac{|\widehat{v} \widehat{v_1} \widehat{v_2}|}{\langle \sigma 
\rangle^{\frac{1}{2}-} \langle \sigma_1 \rangle^{\frac{1}{2}} \langle \sigma_2 
\rangle^{\frac{1}{2}}} $$
by our assumption $k \ge \frac{l+2}{2}$ . This integral is treated by Lemma 
\ref{LemmaB} as before. \\
{\bf Region B:} $|\xi_1| \ge 2 |\xi_2|$ $(\Rightarrow \, |\xi| \sim |\xi_1|)$ 
(and similarly $|\xi_2| \ge 2|\xi_1|$). \\
Using $k \le l+2$ w.l.o.g. and (\ref{3.1}) we get
$$
W \le c \int \frac{|\widehat{v} \widehat{v_1} \widehat{v_2}| \langle \xi_1 
\rangle^{l+2-k}}{\langle \xi_2 \rangle^k  \langle \sigma \rangle^{\frac{1}{2}-} 
\langle \sigma_1 \rangle^{\frac{1}{2}} \langle \sigma_2 \rangle^{\frac{1}{2}}}
\le c \int \frac{|\widehat{v} \widehat{v_1} \widehat{v_2}| (\langle \sigma 
\rangle + \langle \sigma_1 \rangle + \langle \sigma_2 
\rangle)^{\frac{l+2-k}{2}}}{\langle \xi_2 \rangle^k  \langle \sigma 
\rangle^{\frac{1}{2}-} \langle \sigma_1 \rangle^{\frac{1}{2}} \langle \sigma_2 
\rangle^{\frac{1}{2}}} \, .
$$ 
The condition $k>l+1$ is required to produce positive exponents of the $\sigma$ 
- modules in the denominator. Moreover we have $k>0$ so that we can apply Lemma 
\ref{LemmaA} with e.g. $a=\frac{1}{2}-\frac{l+2-k}{2}-$ , $a_1 = a_2 
=\frac{1}{2}$ and $m=k$ , so that $2(a+a_1+a_2)+m=k+(k-l)+1- > \frac{5}{2}$ , 
because $k-l>1$ and $k \ge \frac{l+2}{2} \ge \frac{1}{2}$ . This completes the 
proof of Lemma \ref{LemmaD}.
\begin{Cor}
\label{Corollary to Lemma D}
Under the assumptions of Lemma \ref{LemmaD} we get for $k \ge 1$:
\begin{eqnarray*} 
\lefteqn{\|(-\Delta)^{\frac{1}{2}} (D\bar{\varphi}_1 D 
\varphi_2)\|_{X_{\pm}^{l+1,-\frac{1}{2}+}}} \\
& \le & c 
T^{\Theta} (\|D \varphi_1\|_{X^{1,\frac{1}{2}}} \|D 
\varphi_2\|_{X^{k,\frac{1}{2}}} + \|D \varphi_1\|_{X^{k,\frac{1}{2}}} \|D 
\varphi_2\|_{X^{1,\frac{1}{2}}}) \, . 
\end{eqnarray*} 
\end{Cor} 
{\bf Proof:} Using Lemma \ref{LemmaD} with $k=1$ , $ l=0-$ we get
\begin{equation}
\label{D}
\| D\bar{\varphi}_1 D \varphi_2\|_{X_{\pm}^{2-,-\frac{1}{2}+}}  \le  c 
T^{\Theta} \|D \varphi_1\|_{X^{1,\frac{1}{2}}} \|D 
\varphi_2\|_{X^{1,\frac{1}{2}}} \, , 
\end{equation}  
which gives as in the proof of Corollary \ref{Corollary to Lemma C} for $l \ge 
0-$ :
\begin{eqnarray*}
\lefteqn{\|(-\Delta)^{\frac{1}{2}} (D\bar{\varphi}_1 D 
\varphi_2)\|_{X_{\pm}^{l+1,-\frac{1}{2}+}}  \le  \| D\bar{\varphi}_1 D 
\varphi_2\|_{X_{\pm}^{l+2,-\frac{1}{2}+}}} \\ 
& \le &  c 
T^{\Theta} (\|D \varphi_1\|_{X^{1,\frac{1}{2}}} \|D 
\varphi_2\|_{X^{l+1+,\frac{1}{2}}} + \|D \varphi_1\|_{X^{l+1+,\frac{1}{2}}} \|D 
\varphi_2\|_{X^{1,\frac{1}{2}}}) \\
& \le &  c 
T^{\Theta} (\|D \varphi_1\|_{X^{1,\frac{1}{2}}} \|D 
\varphi_2\|_{X^{k,\frac{1}{2}}} + \|D \varphi_1\|_{X^{k,\frac{1}{2}}} \|D 
\varphi_2\|_{X^{1,\frac{1}{2}}}) \, ,
\end{eqnarray*} 
whereas for $ l \le 0- $ we get obviously by (\ref{D}):
\begin{eqnarray*}
\lefteqn{\hspace{-4em} \|(-\Delta)^{\frac{1}{2}} D\bar{\varphi}_1 D 
\varphi_2\|_{X_{\pm}^{l+1,-\frac{1}{2}+}} \le \| D\bar{\varphi}_1 D 
\varphi_2\|_{X_{\pm}^{2-,-\frac{1}{2}+}}  \le  c 
T^{\Theta} \|D \varphi_1\|_{X^{1,\frac{1}{2}}} \|D 
\varphi_2\|_{X^{1,\frac{1}{2}}}} \\
& \le &  c 
T^{\Theta} (\|D \varphi_1\|_{X^{1,\frac{1}{2}}} \|D 
\varphi_2\|_{X^{k,\frac{1}{2}}} + \|D \varphi_1\|_{X^{k,\frac{1}{2}}} \|D 
\varphi_2\|_{X^{1,\frac{1}{2}}}) \, .
\end{eqnarray*} 
\begin{lemma}
\label{LemmaC'}
Let $n=3$ , $ l \ge -1$ , $ l+1 \le k \le l+2 $ , and let $\varphi$ , $\chi$ be 
given with support in $\{|t| \le cT \}$ . Then the following estimate holds: 
$$ \|(-\Delta)^{-\frac{1}{2}} (D\varphi D\chi)\|_{X^{k,-\frac{1}{2}}} \le 
cT^{\Theta} \|D\varphi\|_{X^{k,\frac{1}{2}}} 
\|D\chi\|_{X_{\pm}^{l,\frac{1}{2}}} $$
with $ \Theta = \Theta(k,l) > 0 $ . 
\end{lemma}
{\bf Remark:} For $l \le 0$ we can obviously replace 
$\|D\chi\|_{X_{\pm}^{l,\frac{1}{2}}}$ by 
$\|\chi\|_{\dot{X}_{\pm}^{l+1,\frac{1}{2}}}$ .\\
{\bf Proof:} We repeat the proof of Lemma \ref{LemmaC} replacing everywhere 
$\langle \sigma_1\rangle^{\frac{1}{2}-}$ by $\langle 
\sigma_1\rangle^{\frac{1}{2}}.$ Then we can allow $(k,l)=(0,-1)$ in Region A, 
Case 1. The strong inequality $k<l+2$ was only used in Region C b. Here the case 
$k=l+2$ is also possible, if  $\langle \sigma_1\rangle^{\frac{1}{2}}$ appears 
instead of $\langle \sigma_1\rangle^{\frac{1}{2}-}$ . Just remark that in the 
limiting case $k=l+2$ we have $k>0$ so that Lemma \ref{LemmaA} can be applied.
\begin{Cor}
\label{Corollary to Lemma C'}
Under the assumptions of Lemma \ref{LemmaC'} we have
$$ \|(-\Delta)^{-\frac{1}{2}} (D\varphi D\chi)\|_{X^{k,-\frac{1}{2}}} \le 
cT^{\Theta} (\|D\varphi\|_{X^{1,\frac{1}{2}}} 
\|\chi\|_{X_{\pm}^{l+1,\frac{1}{2}}} + \|D\varphi\|_{X^{k,\frac{1}{2}}} 
\|\chi\|_{X_{\pm}^{0,\frac{1}{2}}}) \, . $$ 
\end{Cor} 
\begin{lemma}
\label{LemmaD'}
Let $n=3$ , $l \ge -1$ , $ k \ge \frac{l+2}{2}$ , $k=l+1$ and suppose 
$\varphi_1$ and $ \varphi_2$ are supported in $\{|t| \le cT\}$ . Then
$$ \|D\bar{\varphi}_1 D \varphi_2\|_{X_{\pm}^{l+2,-\frac{1}{2}}} \le c 
T^{\Theta} \|D \varphi_1\|_{X^{k,\frac{1}{2}+}} \|D 
\varphi_2\|_{X^{k,\frac{1}{2}+}} $$
with $ \Theta = \Theta(k,l) > 0 $ .
\end{lemma} 
{\bf Remark:} We can replace $X_{\pm}^{l+2,-\frac{1}{2}}$ by 
$\dot{X}_{\pm}^{l+2,-\frac{1}{2}}$ . \\
{\bf Proof:} Replacing $\langle \sigma \rangle^{\frac{1}{2}-}$ by $\langle 
\sigma \rangle^{\frac{1}{2}}$ and $\langle \sigma_i \rangle^{\frac{1}{2}}$
by $\langle \sigma_i \rangle^{\frac{1}{2}+}$ everywhere we repeat the proof of 
Lemma \ref{LemmaD}. The strong condition $k>l+1$ was only required in Region B 
to produce positive exponents of the $\sigma$ - modules in the denominator. In 
the limiting case $k=l+1$ (remark that $k>0$ here) we use Lemma \ref{LemmaA} 
with e.g. $a=0$ , $a_1=\frac{1}{2}+$ , $a_2=\frac{1}{2}+$ and $m=k$ and get the 
inequality
\begin{equation}
\label{****}
2(a+a_1+a_2)+m = 2+k+ > \frac{5}{2} \, ,
\end{equation}
if $k \ge \frac{1}{2}$ . This completes the proof. \\
{\bf Remark:} For $k> \frac{1}{2}$ we can replace $X^{k,\frac{1}{2}+}$ by 
$X^{k,\frac{1}{2}}$ in the statement of Lemma \ref{LemmaD'}. \\
This follows immediately, because in this case condition (\ref{****}) with $a=0$ 
, $ a_1=a_2=\frac{1}{2}$ is also satisfied.
\begin{Cor}
\label{Corollary to Lemma D'}
Under the assumptions of Lemma \ref{LemmaD'} and $k \ge 1$ we get 
\begin{eqnarray*} 
\lefteqn{\|(-\Delta)^{\frac{1}{2}} (D\bar{\varphi}_1 D 
\varphi_2)\|_{X_{\pm}^{l+1,-\frac{1}{2}}}} \\
& \le & c 
T^{\Theta} (\|D \varphi_1\|_{X^{1,\frac{1}{2}}} \|D 
\varphi_2\|_{X^{k,\frac{1}{2}}} + \|D \varphi_1\|_{X^{k,\frac{1}{2}}} \|D 
\varphi_2\|_{X^{1,\frac{1}{2}}}) \, . 
\end{eqnarray*} 
\end{Cor}

Because we were forced to replace $X^{k,-\frac{1}{2}+}$ by $X^{k,-\frac{1}{2}}$ 
in the limiting case $k=l+2$ in Lemma \ref{LemmaC} we have to give an additional 
estimates where $X^{k,-\frac{1}{2}}$ is replaced by $Y^k$ (in order to apply 
(\ref{b}) later). Similarly, because $X_{\pm}^{l+2,-\frac{1}{2}+}$ had to be 
replaced by $X_{\pm}^{l+2,-\frac{1}{2}}$ in the limiting case $k=l+1$ in Lemma 
\ref{LemmaD} we need an estimate where $X_{\pm}^{l+2,-\frac{1}{2}}$ is replaced 
by $Y_{\pm}^{l+2}$ .
\begin{lemma}
\label{LemmaC''}
Let $n=3$ , $ l \ge -1$ , $l+1 \le k \le l+2$ be given and let $\varphi$ and 
$\chi$ be supported in $\{|t| \le cT\}$ . Then 
$$ \|(-\Delta)^{-\frac{1}{2}} (D\varphi D\chi)\|_{Y^k} \le 
cT^{\Theta} \|D\varphi\|_{X^{k,\frac{1}{2}}} 
\|D\chi\|_{X_{\pm}^{l,\frac{1}{2}}} $$
with $ \Theta = \Theta(k,l) > 0 $ . 
\end{lemma}
{\bf Remark:} For $l \le 0$ we can replace $\|D\chi\|_{X_{\pm}^{l,\frac{1}{2}}}$ 
by $\|\chi\|_{\dot{X}_{\pm}^{l+1,\frac{1}{2}}}$ . 
\begin{Cor}
\label{Corollary to Lemma C''}
Under the assumptions of Lemma \ref{LemmaC''} we have
$$ \|(-\Delta)^{-\frac{1}{2}} (D\varphi D\chi)\|_{Y^k} \le 
cT^{\Theta} (\|D\varphi\|_{X^{1,\frac{1}{2}}} 
\|\chi\|_{X_{\pm}^{l+1,\frac{1}{2}}} + \|D\varphi\|_{X^{k,\frac{1}{2}}} 
\|\chi\|_{X_{\pm}^{0,\frac{1}{2}}}) \, . $$ 
\end{Cor}
{\bf Proof of Lemma \ref{LemmaC''}:}
Defining $v$ and $v_2$ as in the proof of Lemma \ref{LemmaC} and 
$\widehat{\psi}(\xi_1) := \langle \xi_1 \rangle^k \widehat{w_1}(\xi_1)$ with 
$w_1 \in L_x^2$ , so that $\psi$ denotes a generic function in $H_x^{-k},$ we 
have to show
$$ \tilde{S}:= \int \frac{|\widehat{v} \widehat{w_1} 
\widehat{v_2}|\,|\xi_1|^{-1} 
\langle \xi_1 \rangle^k}{\langle \xi \rangle^l \langle \xi_2 \rangle^k \langle 
\sigma \rangle^{\frac{1}{2}} \langle \sigma_1 \rangle \langle \sigma_2 
\rangle^{\frac{1}{2}}} \le c T^{\Theta} 
\|v\|_{L^2_{xt}}\|w_1\|_{L^2_{x}}\|v_2\|_{L^2_{xt}} \, . $$
The only case where the strict inequality $k<l+2$ was used in the proof of Lemma 
\ref{LemmaC} was the region $|\xi_1| \ge 2|\xi_2|$ and $|\xi_1| \ge 1$ . In all 
other regions we define $\widehat{v_1} := \langle \sigma_1 
\rangle^{-\frac{1}{2}-} \widehat{w_1}$. Then one easily checks 
$\|v_1\|_{L^2_{xt}} \le c \|w_1\|_{L_x^2}$ and $\tilde{S}$ can be replaced by
$$ \int \frac{|\widehat{v} \widehat{v_1} \widehat{v_2}|\,|\xi_1|^{-1} \langle 
\xi_1 \rangle^k}{\langle \xi \rangle^l \langle \xi_2 \rangle^k \langle \sigma 
\rangle^{\frac{1}{2}} \langle \sigma_1 \rangle^{\frac{1}{2}-} \langle \sigma_2 
\rangle^{\frac{1}{2}}} \, . $$
This is exactly the integral treated in the proof of Lemma \ref{LemmaC}, so that 
the desired result in these regions follows using the remarks to Lemma 
\ref{LemmaA} and Lemma \ref{LemmaB} taking into account that $w_1$ fulfills no 
support property. \\
It remains to consider the region where $|\xi_1| \ge 2|\xi_2|$ and $|\xi_1| \ge 
1$ and $l+1 \le k \le l+2$ . In this case we get as in Lemma \ref{LemmaC}
\begin{eqnarray*}
\tilde{S} & \le & c \int \frac{|\widehat{v} \widehat{w_1} \widehat{v_2}| \langle 
\xi_1 \rangle^{k-l-1}}{\langle \xi_2 \rangle^k \langle \sigma 
\rangle^{\frac{1}{2}} \langle \sigma_1 \rangle \langle \sigma_2 
\rangle^{\frac{1}{2}}} \\
& \le & c \int \frac{|\widehat{v} \widehat{w_1} \widehat{v_2}| (\langle \sigma 
\rangle + \langle \sigma_2 \rangle + \langle \sigma_1 \rangle 
\phi_{\{c_1|\sigma_1| \le |\xi_1|^2 \le 
c_2|\sigma_1|\}})^{\frac{k-l-1}{2}}}{\langle \xi_2 \rangle^k \langle \sigma 
\rangle^{\frac{1}{2}} \langle \sigma_1 \rangle \langle \sigma_2 
\rangle^{\frac{1}{2}}} \, . 
\end{eqnarray*} 
Here we used (\ref{3.2}). The two terms coming from $\langle \sigma \rangle$ and 
$\langle \sigma_2 \rangle$ in the numerator are treated by defining 
$\widehat{v_1}$ as before by Lemma \ref{LemmaA} with e.g. 
$a=\frac{1}{2}-\frac{k-l-1}{2} \ge 0$ , $a_1=\frac{1}{2}-$ , $a_2=\frac{1}{2}$ , 
$ m=k$ , which implies $2(a+a_1+a_2)+m=4+l- > \frac{5}{2}$ , whereas the term 
coming from $\langle \sigma_1 \rangle$ is treated by defining $\widehat{v_1} 
:=\langle \sigma_1 \rangle^{-\frac{1}{2}} \widehat{w_1} 
\phi_{\{c_1|\sigma_1| \le |\xi_1|^2 \le c_2|\sigma_1|\}}.$ One can easily 
show $\|v_1\|_{L^2_{xt}} \le c \|w_1\|_{L^2_x}$ , so that we only have to give 
the estimate
$$\int \frac{|\widehat{v} \widehat{v_1} \widehat{v_2}|}{\langle \xi_2 \rangle^k 
\langle \sigma \rangle^{\frac{1}{2}} \langle \sigma_1 
\rangle^{\frac{1}{2}-\frac{k-l-1}{2}} \langle \sigma_2 \rangle^{\frac{1}{2}}} 
\le c T^{\Theta} \|v\|_{L^2_{xt}} \|v_1\|_{L^2_{xt}} \|v_2\|_{L^2_{xt}} \, . $$
This can be done by Lemma \ref{LemmaA} (+ remark) with $a=a_2=\frac{1}{2}$ , 
$a_1=\frac{1}{2}-\frac{k-l-1}{2} \ge 0$ and $m=k$ which implies $2(a+a_1+a_2)+m 
= 4-l > \frac{5}{2}$ .

We also get
\begin{lemma}
\label{LemmaD''}
Let $n=3$ , $l \ge -1$ , $ k \ge \frac{l+2}{2}$ , $k=l+1$ and suppose 
$\varphi_1$ and $ \varphi_2$ are supported in $\{|t| \le cT\}$ . Then
$$ \|D\bar{\varphi}_1 D \varphi_2\|_{Y_{\pm}^{l+2}} \le c 
T^{\Theta} \|D \varphi_1\|_{X^{k,\frac{1}{2}+}} \|D 
\varphi_2\|_{X^{k,\frac{1}{2}+}} $$
with $ \Theta = \Theta(k,l) > 0 $ . If $k > \frac{1}{2}$ , we can replace 
$X^{k,\frac{1}{2}+}$ by $X^{k,\frac{1}{2}}$ . 
\end{lemma} 
{\bf Remark:} We can obviously replace $Y_{\pm}^{l+2}$ by $\dot{Y}_{\pm}^{l+2}$ 
and $k=l+1$ by $k \ge l+1$.
{\bf Proof:} Defining $v_1$ and $v_2$ similarly as in the proof of Lemma 
\ref{LemmaD} and $\widehat{\psi}(\xi):= \langle \xi \rangle^{l+2} 
\widehat{w}(\xi)$ 
with $w \in L_x^2$ (so that $\psi$ is a generic function in $H_x^{-l-2}$), we 
have to show for any $\epsilon > 0$:
$$  \tilde{W} := \int \frac{|\widehat{w} \widehat{v_1} \widehat{v_2}| \langle 
\xi 
\rangle^{l+2}}{\langle \xi_1 \rangle^k \langle \xi_2 \rangle^k  \langle \sigma 
\rangle \langle \sigma_1 \rangle^{\frac{1}{2}+\epsilon} \langle \sigma_2 
\rangle^{\frac{1}{2}+\epsilon}} \le c T^{\Theta} \|w\|_{L^2_{x}} 
\|v_1\|_{L^2_{xt}} \|v_2\|_{L^2_{xt}} \, . $$ 
In region A of the proof of Lemma \ref{LemmaD} we define $\widehat{v}:=\langle 
\sigma \rangle^{-\frac{1}{2}-\frac{\epsilon}{2}} \widehat{w}$ such that 
$\|v\|_{L^2_{xt}} \le c \|w\|_{L_x^2}$ and $\tilde{W}$ is estimated by
$$ c\int \frac{|\widehat{v} \widehat{v_1} \widehat{v_2}| \langle \xi_1 
\rangle^{l+2-2k}}{\langle \sigma 
\rangle^{\frac{1}{2}-\frac{\epsilon}{2}} \langle \sigma_1 
\rangle^{\frac{1}{2}+\epsilon} \langle \sigma_2 
\rangle^{\frac{1}{2}+\epsilon}} \le c\int \frac{|\widehat{v} \widehat{v_1} 
\widehat{v_2}|}{\langle \sigma 
\rangle^{\frac{1}{2}-\frac{\epsilon}{2}} \langle \sigma_1 
\rangle^{\frac{1}{2}+\epsilon} \langle \sigma_2 
\rangle^{\frac{1}{2}+\epsilon}} \, , $$
which can be estimated by $cT^{\Theta} \|v\|_{L^2} \|v_1\|_{L^2} \|v_2\|_{L^2} $ 
by Lemma \ref{LemmaB} (+ remark) as before. 
In region B of the proof of Lemma \ref{LemmaD} we get using $k=l+1$ and 
(\ref{3.3}) :
\begin{eqnarray*}
\tilde{W} & \le & c \int \frac{|\widehat{w} \widehat{v_1} \widehat{v_2}| \langle 
\xi_1 \rangle}{\langle \xi_2 \rangle^k  \langle \sigma 
\rangle \langle \sigma_1 \rangle^{\frac{1}{2}+\epsilon} \langle \sigma_2 
\rangle^{\frac{1}{2}+\epsilon}} \\
& \le & c \int \frac{|\widehat{w} \widehat{v_1} \widehat{v_2}| (\langle \sigma_1 
\rangle + \langle \sigma_2 \rangle + \langle \sigma \rangle
 \phi_{\{ c_1 |\sigma| \le |\xi|^2 \le c_2 |\sigma| \}}
)^{\frac{1}{2}}}{\langle \xi_2 \rangle^k  \langle \sigma 
\rangle \langle \sigma_1 \rangle^{\frac{1}{2}+\epsilon} \langle \sigma_2 
\rangle^{\frac{1}{2}+\epsilon}} \, .
\end{eqnarray*}
The two terms coming from $\langle \sigma_1 \rangle$ and $\langle \sigma_2 
\rangle $ in the numerator are treated by defining $\widehat{v}$ as before by 
Lemma \ref{LemmaA} with e.g. $a_1 = \epsilon$ , $a_2 =\frac{1}{2} + \epsilon$ , 
$a=\frac{1}{2}-\frac{\epsilon}{2}$ , $ m=k \ge \frac{1}{2},$ so that 
\begin{equation}
\label{*****}
2(a+a_1+a_2)+m > \frac{5}{2} \, . 
\end{equation}
The term coming from $\langle \sigma \rangle$ is treated by defining 
$\widehat{v} := \langle \sigma 
\rangle^{-\frac{1}{2}} \widehat{w} \phi_{ \{c_1|\sigma| \le |\xi|^2 \le 
c_2|\sigma| \} }$ , so that $ \|v\|_{L^2_{xt}} \le c \|w\|_{L_x^2} $ . 
Thus it remains to show
$$ \int \frac{|\widehat{v} \widehat{v_1} \widehat{v_2}|}{\langle \xi_2 \rangle^k 
\langle \sigma_1 \rangle^{\frac{1}{2}+\epsilon} \langle \sigma_2 
\rangle^{\frac{1}{2}+\epsilon}} \le cT^{\Theta} \|v\|_{L^2} \|v_1\|_{L^2} 
\|v_2\|_{L^2} \, . $$
This is true by Lemma \ref{LemmaA} with $a_1=a_2=\frac{1}{2}+\epsilon$ , $a=0$ , 
$m=k \ge \frac{1}{2}$ , thus $2(a+a_1+a_2)+m > \frac{5}{2}$ . \\
If $k > \frac{1}{2}$ , we can easily modify the proof by replacing $\langle 
\sigma_j \rangle^{\frac{1}{2}+\epsilon}$ by $\langle 
\sigma_j\rangle^{\frac{1}{2}}$ $(j=1,2),$ because the decisive condition 
(\ref{*****}) in this case also holds. 
\begin{Cor}
\label{Corollary to Lemma D''}
Under the assumptions of Lemma \ref{LemmaD''} and $k \ge 1$ we get 
\begin{eqnarray*} 
\lefteqn{\|(-\Delta)^{\frac{1}{2}} (D\bar{\varphi}_1 D 
\varphi_2)\|_{Y_{\pm}^{l+1}}} \\
& \le & c 
T^{\Theta} (\|D \varphi_1\|_{X^{1,\frac{1}{2}}} \|D 
\varphi_2\|_{X^{k,\frac{1}{2}}} + \|D \varphi_1\|_{X^{k,\frac{1}{2}}} \|D 
\varphi_2\|_{X^{1,\frac{1}{2}}}) \, . 
\end{eqnarray*} 
\end{Cor} 
{\bf Proof:} follows from Lemma \ref{LemmaD''} and the remark to that Lemma. 
\begin{theorem}
\label{Theorem1}
In space dimension $n=3$ assume $l \ge -1$ , $l+1 \le k \le l+2, $  $ k \ge 
\frac{l+2}{2}$ , and
$$B\varphi_0 \in H^k({\bf R^3}) \, , \, B\chi_0 \in H^l({\bf R}^3) \, , \, 
\chi_1 \in H^l({\bf R}^3) \, . $$
Then there exists $1 \ge T > 0$ , 
$T=T(\|B\varphi_0\|_{H^k},\|B\chi_0\|_{H^l},\|\chi_1\|_{H^l})$ , such that the 
problem (\ref{1.1}),(\ref{1.2}),(\ref{1.3}) has a unique solution 
$(\varphi,\chi)$ with 
$$B\varphi \in X^{k,b}[0,T] \quad , \quad B\chi \, , \, \chi_t \in 
X_+^{l,b_1}[0,T] + X_-^{l,b_1}[0,T] \, .$$ 
Here $b=\frac{1}{2}+$ , $b_1=\frac{1}{2}+$ , if $l+1<k<l+2$ , $b=\frac{1}{2}$ , 
$b_1=\frac{1}{2}+$ , if $ k=l+2$ , and $b=\frac{1}{2}+$ , $b_1=\frac{1}{2}$ , if 
$k=l+1$ . This solution satisfies 
$$B\varphi \in C^0([0,T],H^k({\bf R}^3)) \quad , \quad B\chi \, , \, \chi_t \in 
C^0([0,T],H^l({\bf R}^3)) \, . $$ 
If $l \le 0$ we can replace $B\chi_0 \, , \, \chi_1 \in H^l$ by $\chi_0 \in 
\dot{H}^{l+1} $ , $\chi_1 \in \dot{H}^l$ ,  and $ B\chi \, , \, \chi_t \in 
X_+^{l,b_1}[0,T] + X_-^{l,b_1}[0,T] $ by $\chi \in \dot{X}_+^{l+1,b_1}[0,T] + 
\dot{X}_-^{l+1,b_1}[0,T]$ , $ \chi_t \in \dot{X}_+^{l,b_1}[0,T] + 
\dot{X}_-^{l,b_1}[0,T]$ , and we have $\chi \in C^0([0,T],\dot{H}^{l+1}({\bf 
R}^3)$ , $\chi_t \in C^0([0,T],\dot{H}^{l}({\bf R}^3)$ .
\end{theorem}     
{\bf Proof:} We replace our system of integral equations by the cut-off system
\begin{eqnarray*}
B\varphi(t) & = & \psi_1(t) B e^{it\Delta} \varphi_0 - \frac{1}{2i} \psi_T(t) 
\int_0^t e^{i(t-s)\Delta} B^{-1}(((\psi_{2T}(s) \nabla \varphi(s)) \\
& & \hspace{12em} \times (\psi_{2T}(s) \nabla(\chi_+(s)+\chi_-(s))) \cdot e) \, 
ds \\
B\chi_{\pm}(t) & = & \psi_1(t) B e^{\pm itB} \chi_{\pm 0} \mp \frac{1}{i} 
\psi_T(t) \int_0^t e^{\mp i(t-s)B} B^2(((\psi_{2T}(s) \nabla \bar{\varphi}(s)) 
\\
& & \hspace{15em} \times (\psi_{2T}(s) \nabla \varphi(s))) \cdot e) \, ds \, ,
\end{eqnarray*}  
which we want to solve globally in $t$. This gives a solution of the original 
system in $[0,T]$ . The factors $\psi_{2T}$ here allow to assume that the 
factors in the nonlinearities are supported in $\{|t| \le 2T\}$ . We want to use 
the contraction mapping principle and consider the case $l+1<k<l+2$ first. \\
The linear parts are treated as follows:
$$ \|\psi_1(t) B e^{it\Delta} \varphi_0\|_{X^{k,b}} \le c \|B\varphi_0\|_{H^k} 
$$
and
$$ \|\psi_1(t) B e^{\pm itB} \chi_{\pm 0}\|_{X_{\pm}^{l,b_1}} \le c \|B\chi_{\pm 
0}\|_{H^l} \, . $$
Using (\ref{a}) the integral term in the first equation can be estimated in the 
$X^{k,\frac{1}{2}+}$ - norm by
$$ cT^{0+} \|B^{-1}(((\psi_{2T} \nabla \varphi) \times (\psi_{2T} \nabla(\chi_+ 
+ \chi_-))) \cdot e)\|_{X^{k,-\frac{1}{2}++}} \, ,  $$ 
which by Lemma \ref{LemmaC} and (\ref{c}) is majorized by
\begin{eqnarray*}
\lefteqn{ cT^{\Theta +} \|B(\psi_{2T} 
\varphi)\|_{X^{k,\frac{1}{2}}}(\|B(\psi_{2T} \chi_+)\|_{X_+^{l,\frac{1}{2}}} + 
\|B(\psi_{2T} \chi_-)\|_{X_-^{l,\frac{1}{2}}}) } \\
& \le & cT^{\Theta -} \|B \varphi\|_{X^{k,\frac{1}{2}}}(\|B 
\chi_+\|_{X_+^{l,\frac{1}{2}}} + \|B \chi_-\|_{X_-^{l,\frac{1}{2}}}) \, ,
\end{eqnarray*}
where $\Theta > 0 $ . \\
The integral term in the second equation can be estimated in the 
$X_{\pm}^{l,\frac{1}{2}+}$ - norm similarly by use of Lemma \ref{LemmaD} instead 
of Lemma \ref{LemmaC} and leads to the bound $ cT^{\Theta -} 
\|B\varphi\|_{X^{k,\frac{1}{2}}}^2$ . \\
The standard contraction argument then gives a unique solution $B\varphi \in 
X^{k,b}$ , $B\chi_{\pm} \in X_{\pm}^{l,b_1}$ of the cut-off system for small 
enough $T$. \\
If $k=l+1$ the estimates for the first equation remain unchanged whereas Lemma 
\ref{LemmaD} is no longer true and forces us to choose $b_1 = \frac{1}{2}$ , so 
that the integral term in the $X_{\pm}^{l,b_1}$ - norm is estimated by (\ref{b}) 
by
$$ \|B^2(((\psi_{2T} \nabla \bar{\varphi}) \times (\psi_{2T} \nabla \varphi)) 
\cdot e))\|_{X_{\pm}^{l,-\frac{1}{2}}} + \|B^2(((\psi_{2T} \nabla \bar{\varphi}) 
\times (\psi_{2T} \nabla \varphi)) \cdot e))\|_{Y^l} \, . $$
The first term can be treated by Lemma \ref{LemmaD'} and (\ref{c}) and gives the 
bound $cT^{\Theta} \|B(\psi_{2T}\varphi)\|_{X^{k,\frac{1}{2}+}}^2 \le cT^{\Theta 
-} \|B\varphi\|_{X^{k,\frac{1}{2}+}}^2 $ , whereas the second term gives the 
same bound by Lemma \ref{LemmaD''}. So we get a unique solution $B\varphi \in 
X^{k,\frac{1}{2}+} $ , $ B\chi_{\pm} \in X_{\pm}^{l,\frac{1}{2}}$ . \\
If $k=l+2$ the estimates for the second equation remain unchanged, whereas Lemma 
\ref{LemmaC} is no longer true and thus requires $b=\frac{1}{2}$ so that the 
integral term in the $X^{k,b}$ - norm is bounded by
\begin{eqnarray*}
\lefteqn{\|B(((\psi_{2T} \nabla \varphi) \times (\psi_{2T} \nabla(\chi_+ + 
\chi_-))) \cdot e)\|_{X^{k,-\frac{1}{2}}}} \\
& + & \|B(((\psi_{2T} \nabla \varphi) \times (\psi_{2T} \nabla(\chi_+ + 
\chi_-))) \cdot e)\|_{Y^k} \, .
\end{eqnarray*}
These terms are treated by Lemma \ref{LemmaC'} and Lemma \ref{LemmaC''}, which 
gives the bound
\begin{eqnarray*}
\lefteqn{ cT^{\Theta} \|B\psi_{2T} \varphi\|_{X^{k,\frac{1}{2}}} (\|B\psi_{2T} 
\chi_+\|_{X^{l_+,\frac{1}{2}}} + \|B\psi_{2T} \chi_-\|_{X^{l_-,\frac{1}{2}}}) } 
\\
& \le & cT^{\Theta -} \|B \varphi\|_{X^{k,\frac{1}{2}}} (\|B 
\chi_+\|_{X^{l_+,\frac{1}{2}}} + \|B \chi_-\|_{X^{l_-,\frac{1}{2}}}) \, ,
\end{eqnarray*}
which leads to a unique solution $B\varphi \in X^{k,\frac{1}{2}}$ , $B\chi_{\pm} 
\in X_{\pm}^{l,\frac{1}{2}+}$ of the cut-off system. \\
To prove uniqueness for the original system of integral equations in $[0,T]$ 
(without cut-offs) let $(\varphi,\chi_{\pm})$ be any solution with $B\varphi \in 
X^{k,b}[0,T]$ , $ B\chi_{\pm} \in X_{\pm}^{l,b_1}[0,T]$ . Consider e.g. the case 
$l+1<k<l+2$ 
and $b=\frac{1}{2}+$ , $b_1=\frac{1}{2}+$ . Let 
$(\tilde{\varphi},\tilde{\chi}_{\pm})$ be any extension with $B\tilde{\varphi} 
\in X^{k,b}$ , $B\tilde{\chi}_{\pm} \in X_{\pm}^{l,b_1}$ . Then we have by the 
same estimates as above:  
\begin{eqnarray*}
\lefteqn{\|\int_0^t e^{i(t-s)\Delta} B^{-1} (((\nabla \varphi(s) \times 
\nabla(\chi_+(s) + \chi_-(s))) \cdot e) \, ds \|_{X^{k,b}[0,T]}}
 \\
&\hspace{-0.6em} \le & \hspace{-0.6em}\|\psi_T(t)\int_0^t e^{i(t-s)\Delta} 
B^{-1} (((\psi_{2T}(s)\nabla 
\tilde{\varphi}(s) \times \psi_{2T}(s)\nabla(\tilde{\chi}_+(s) + 
\tilde{\chi}_-(s))) 
\cdot e) \, ds \|_{X^{k,b}} \\
& \hspace{-0.6em}\le & \hspace{-0.6em}cT^{\Theta -} 
\|B\tilde{\varphi}\|_{X^{k,\frac{1}{2}}} 
(\|B\tilde{\chi}_+\|_{X^{l,\frac{1}{2}}_+} + 
\|B\tilde{\chi}_-\|_{X^{l,\frac{1}{2}}_-}) \, . 
\end{eqnarray*} 
Thus 
\begin{eqnarray*}
\lefteqn{\|\int_0^t e^{i(t-s)\Delta} B^{-1} (((\nabla \varphi(s) \times 
\nabla(\chi_+(s) + \chi_-(s))) \cdot e) \, ds \|_{X^{k,b}[0,T]}}
 \\
& \le & cT^{\Theta -} \|B\varphi\|_{X^{k,\frac{1}{2}}[0,T]} 
(\|B\chi_+\|_{X^{l,\frac{1}{2}}_+[0,T]} + 
\|B\chi_-\|_{X^{l,\frac{1}{2}}_-[0,T]}) \, . 
\end{eqnarray*} 
Similarly we can treat this term in the other cases using the $Y$ - spaces and 
also the integral term in the second integral equation. A standard argument 
implies uniqueness for the original system in $[0,T]$ .\\
The claim that $B\varphi$ belongs to $C^0([0,T],H^k)$ and $B\chi$ to 
$C^0([0,t],H^l)$ follows directly from the embeddings $\, X^{k,b}[0,T] \subset 
C^0([0,T],H^k) \,\,$ and $\, \,X_{\pm}^{l,b_1}[0,T] \subset $ $ C^0([0,T], H^l)$ 
for $b>\frac{1}{2}$ and $b_1 > \frac{1}{2}$ . If $b=\frac{1}{2}$ (or similarly 
$b_1 = \frac{1}{2}$) this follows from the fact that the nonlinearity 
$B^{-1}(((\psi_{2T}\nabla \varphi) \times (\psi_{2T} \nabla \chi_{\pm})) \cdot 
e) $ belongs to $Y^k$ for $B\varphi \in X^{k,\frac{1}{2}}$ and $B\chi_{\pm} \in 
X_{\pm}^{l,\frac{1}{2}}$ (cf. estimate above). This implies by \cite{GTV}, Lemma 
2.2: $\int_0^t e^{i(t-s)\Delta} B(((\psi_{2T} \nabla \varphi) \times (\psi_{2T} 
\nabla (\chi_+ + \chi_-))) \cdot e) \, ds \in C^0({\bf R},H^k({\bf R}^3)), $  
which by the integral equation implies $B\varphi \in C^0([0,T],H^k({\bf R}^3))$ 
. \\
The additional claim for $l \le 0$ follows easily by replacing in the 
application of Lemma \ref{LemmaC}, Lemma \ref{LemmaC'} and Lemma \ref{LemmaC''} 
$\|B\chi_{\pm}\|_{X_{\pm}^{l,\frac{1}{2}}}$ by 
$\|\chi_{\pm}\|_{\dot{X}_{\pm}^{l+1,\frac{1}{2}}}$ and in the application of 
Lemma \ref{LemmaD} and Lemma \ref{LemmaD'} $\|D\bar{\varphi} 
D\varphi\|_{X_{\pm}^{l+2,-\frac{1}{2}(+)}}$ by $\|D\bar{\varphi} 
D\varphi\|_{\dot{X}_{\pm}^{l+2,-\frac{1}{2}(+)}}$ and in Lemma \ref{LemmaD''} 
$\|D\bar{\varphi} D\varphi\|_{Y_{\pm}^{l+2}}$ by $\|D\bar{\varphi} 
D\varphi\|_{\dot{Y}_{\pm}^{l+2}}$ . \\
{\bf Remark:} The case $k=1$ , $l=-1$ especially shows that, given data 
$\varphi_0$ , $ \chi_0$ with $B\varphi_0 \in H^1({\bf R}^3)$ and $\chi_0$ , 
$B^{-1}\chi_1 \in L^2({\bf R}^3)$ , there exists a unique local solution 
$(\varphi,\chi)$ of problem (\ref{1.1}),(\ref{1.2}),(\ref{1.3}) on $[0,T]$ , 
$T=T(\|B\varphi_0\|_{H^1},\|\chi_0\|_{L^2},\|B^{-1}\chi_1\|_{L^2}),$ with 
$B\varphi \in X^{1,\frac{1}{2}}[0,T] $ and $ \chi$ , $B^{-1} \chi_t \in  
X_+^{0,\frac{1}{2}+}[0,T] + X_-^{0,\frac{1}{2}+}[0,T].$ Moreover $B\varphi \in 
C^0([0,T],H^1({\bf R}^3)) $ and $ \chi,B^{-1}\chi_t \in C^0([0,T],L^2({\bf 
R}^3))$ .

Combining the last remark with Proposition \ref{Proposition} we immediately get
\begin{theorem}
\label{Theorem2}
Let $\varphi$ , $\chi_0$ , $\chi_1$ be given with
$$ \|B\varphi_0\|_{H^1} + \|\chi_0\|_{L^2} + \|B^{-1} \chi_1\|_{L^2} < 
\epsilon_0 \, $$
where $\epsilon_0$ is a sufficiently small constant (depending only on $e \in 
{\bf R}^3$ and a Sobolev embedding constant). Then the Cauchy problem 
(\ref{1.1}),(\ref{1.2}),(\ref{1.3}) has a unique global solution 
$(\varphi,\chi)$ with 
$$B\varphi \in X^{1,\frac{1}{2}} \quad , \quad \chi , B^{-1} \chi_t \in  
X_+^{0,\frac{1}{2}+} + X_-^{0,\frac{1}{2}+} \, . $$
Moreover 
$$B\varphi \in C^0({\bf R},H^1({\bf R}^3)) \quad , \quad \chi,B^{-1}\chi_t \in 
C^0({\bf R},L^2({\bf R}^3)) \, .$$ 
\end{theorem}

Using the refinements of the nonlinear estimates given in Corollary 
\ref{Corollary to Lemma C}, Corollary \ref{Corollary to Lemma D}, Corollary 
\ref{Corollary to Lemma C'}, Corollary \ref{Corollary to Lemma D'}, Corollary 
\ref{Corollary to Lemma C''} and Corollary \ref{Corollary to Lemma D''} we get 
the following variant of Theorem \ref{Theorem2}.
\begin{theorem}
\label{Theorem3}
Assume $k \ge 1$ , $ l \ge -1$ , $l+1 \le k \le l+2 $ and 
$$B\varphi_0 \in H^k({\bf R}^3) \quad , \quad \chi_0 , B^{-1}\chi_1 \in 
H^{l+1}({\bf R}^3) \, . $$
Then there exists $1 \ge T > 0 $ , 
$T=T(\|B\varphi_0\|_{H^1},\|\chi_0\|_{L^2},\|B^{-1}\chi_1\|_{L^2})$ ,  such that 
problem (\ref{1.1}), (\ref{1.2}), (\ref{1.3}) has a unique solution 
$(\varphi,\chi)$ with
$$B\varphi \in X^{k,\frac{1}{2}}[0,T] \quad , \quad \chi,B^{-1}\chi_t \in 
X_+^{l+1,b_1}[0,T] + X_-^{l+1,b_1}[0,T] \, ,$$ 
where $b_1=\frac{1}{2}+$ , if $l+1 < k \le l+2$ , and $b_1 =\frac{1}{2}$ , if 
$k=l+1$ . This solution satisfies 
$$B\varphi \in C^0([0,T],H^k({\bf R}^3)) \quad ,  \quad B\chi , B^{-1} \chi_t 
\in C^0([0,T],H^{l+1}({\bf R}^3)) \, .$$ 
\end{theorem}
{\bf Proof:} One has to modify the usual contraction argument in the proof of 
Theorem \ref{Theorem1} combining the following fundamental estimates, which e.g. 
in the case $l+1<k\le l+2$ read as follows:
\begin{eqnarray}
\label{A*}
\|B^{-1}(D\varphi D\chi)\|_{X^{1,-\frac{1}{2}}}  & \hspace{-0.6em} \le & 
\hspace{-0.6em} c T^{\Theta} \|D\varphi\|_{X^{1,\frac{1}{2}}} 
\|\chi\|_{X_{\pm}^{0,\frac{1}{2}}} \\
\label{B*}
\|B^{-1}(D\varphi D\chi)\|_{Y^1} & \hspace{-0.6em} \le & \hspace{-0.6em} c 
T^{\Theta} \|D\varphi\|_{X^{1,\frac{1}{2}}} \|\chi\|_{X_{\pm}^{0,\frac{1}{2}}} 
\\
\label{C*}
\|B^{-1}(D\varphi D\chi)\|_{X^{k,-\frac{1}{2}}} &\hspace{-0.6em} \le & 
\hspace{-0.6em} c T^{\Theta} (\|D\varphi\|_{X^{1,\frac{1}{2}}} 
\|\chi\|_{X_{\pm}^{l+1,\frac{1}{2}}} + \|D\varphi\|_{X^{k,\frac{1}{2}}} 
\|\chi\|_{X_{\pm}^{0,\frac{1}{2}}})\\
\label{D*}
\|B(D\bar{\varphi}_1 D\varphi_2)\|_{X_{\pm}^{1,-\frac{1}{2}+}} & \hspace{-0.6em} 
\le & \hspace{-0.6em} c T^{\Theta} \|D\varphi_1\|_{X^{1,\frac{1}{2}}} 
\|D\varphi_2\|_{X^{1,\frac{1}{2}}}  \\
\nonumber
\|B(D\bar{\varphi}_1 D\varphi_2)\|_{X_{\pm}^{l+1,-\frac{1}{2}+}} & 
\hspace{-0.6em} \le & \hspace{-0.6em} c T^{\Theta} 
(\|D\varphi_1\|_{X^{1,\frac{1}{2}}} \|D\varphi_2\|_{X^{k,\frac{1}{2}}}  +  
\|D\varphi_1\|_{X^{k,\frac{1}{2}}} \|D\varphi_2\|_{X^{1,\frac{1}{2}}})  . \\
\label{E*}
& & 
\end{eqnarray}
Here (\ref{A*}),(\ref{B*}),(\ref{C*}),(\ref{D*}) and (\ref{E*}) follow from 
Lemma \ref{LemmaC'} (+ remark), Lemma \ref{LemmaC''} (+ remark), Corollary 
\ref{Corollary to Lemma C'}, Lemma \ref{LemmaD} and Corollary \ref{Corollary to 
Lemma D}, respectively. \\
In the limiting case $k=l+1$ we only have to replace (\ref{E*}) by
$$ \|B(D\bar{\varphi}_1 D\varphi_2)\|_{X_{\pm}^{l+1,-\frac{1}{2}}}  \le  c 
T^{\Theta} (\|D\varphi_1\|_{X^{1,\frac{1}{2}}} 
\|D\varphi_2\|_{X^{k,\frac{1}{2}}} + \|D\varphi_1\|_{X^{k,\frac{1}{2}}} 
\|D\varphi_2\|_{X^{1,\frac{1}{2}}}) \, , $$
which follows from Corollary \ref{Corollary to Lemma D'}, and to add 
$$ \|B(D\bar{\varphi}_1 D\varphi_2)\|_{Y_{\pm}^{l+1}}  \le  c T^{\Theta} 
(\|D\varphi_1\|_{X^{1,\frac{1}{2}}} \|D\varphi_2\|_{X^{k,\frac{1}{2}}} + 
\|D\varphi_1\|_{X^{k,\frac{1}{2}}} \|D\varphi_2\|_{X^{1,\frac{1}{2}}}) \, , $$
coming from Corollary \ref{Corollary to Lemma D''}. \\
We omit the proof and just refer to \cite{G}, Theorem 1.1, where a detailed 
proof can be found.

Combining Theorem \ref{Theorem3} with Proposition \ref{Proposition} we can also 
show global well-posedness for smoother data, namely
\begin{theorem}
\label{Theorem4}
Assume $k \ge 1$ , $ l \ge -1$ , $l+1 \le k \le l+2 $ and 
$$B\varphi_0 \in H^k({\bf R}^3)\quad , \quad \chi_0,B^{-1}\chi_1 \in 
H^{l+1}({\bf R}^3)$$
 with
$$ \|B\varphi_0\|_{H^1} + \|\chi_0\|_{L^2} + \|B^{-1}\chi_1\|_{L^2} < \epsilon_0 
\, , $$
where $\epsilon_0$ is sufficiently small, dependent only on $e \in {\bf R}^3$ 
and a Sobolev embedding constant. Then the Cauchy problem 
(\ref{1.1}),(\ref{1.2}),(\ref{1.3}) has a unique global solution 
$(\varphi,\chi)$ with 
$$B\varphi \in X^{k,\frac{1}{2}} \quad , \quad \chi,B^{-1}\chi_t \in 
X_+^{l+1,b_1} + X_-^{l+1,b_1} \, , $$
where $b_1=\frac{1}{2}+$ , if $l+1 < k \le l+2$ , and $b_1 = \frac{1}{2}$ , if 
$k=l+1$ . This solution satisfies 
$$B\varphi \in C^0({\bf R},H^k({\bf R}^3)) \quad , \quad \chi,B^{-1}\chi_t \in 
C^0({\bf R},H^{l+1}({\bf R}^3)) \, . $$
\end{theorem}
\section{Local existence in 2+1 dimensions}
\begin{lemma}
\label{LemmaE}
In space dimension $n=2$ the following estimate holds under the assumptions of 
Lemma \ref{LemmaC}:
$$ \|B^{-1+\epsilon}(D\varphi D\chi)\|_{X^{k-\epsilon,-\frac{1}{2}+}} \le c 
T^{\Theta} \|B^{\epsilon} D\varphi\|_{X^{k-\epsilon,\frac{1}{2}}} \|B^{-\delta} 
D\chi\|_{X_{\pm}^{l+\delta,\frac{1}{2}}} $$
with $\Theta > 0$ , if $0 < \epsilon < 1$ and $\delta > 0$ .
\end{lemma}
{\bf Remark:} If $l<0$ , we can replace $\|B^{-\delta} 
D\chi\|_{X_{\pm}^{l+\delta,\frac{1}{2}}}$ by 
$\|\chi\|_{\dot{X}_{\pm}^{l+1,\frac{1}{2}}} $ . \\
{\bf Proof:} We follow the proof of Lemma \ref{LemmaC} and have to give the 
estimate
$$
S:=\left| \int \frac{\widehat{v} \widehat{v_1} \widehat{v_2} 
|\xi_1|^{-1+\epsilon} 
\langle \xi_1 \rangle^{k-\epsilon}}{|\xi|^{-\delta}\langle \xi 
\rangle^{l+\delta} |\xi_2|^{\epsilon} \langle \xi_2 \rangle^{k-\epsilon} \langle 
\sigma \rangle^{\frac{1}{2}} \langle \sigma_1 \rangle^{\frac{1}{2}-} \langle 
\sigma_2 \rangle^{\frac{1}{2}}} \right| \le cT^{\Theta} \|v\|_{L^2_{xt}} 
\|v_1\|_{L^2_{xt}} \|v_2\|_{L^2_{xt}} \, .
$$  
{\bf Region A:} $|\xi_1| \le \frac{1}{2} |\xi_2|$ $\, (\Longrightarrow |\xi| 
\sim|\xi_2|)$ . \\
{\bf Case 1:} $|\xi_1| \ge 1$ , $ |\xi_2| \ge 1$ . \\
The same calculation as in Lemma \ref{LemmaC} gives the desired estimate.\\
{\bf Case 2:} $|\xi_1| \le 1$ , $ |\xi_2| \ge 1 $ . \\
We have
$$ S \le c\int \frac{|\widehat{v} \widehat{v_1} \widehat{v_2}| 
|\xi_1|^{-1+\epsilon}}{\langle \xi_2 \rangle^{l+k} \langle 
\sigma \rangle^{\frac{1}{2}} \langle \sigma_1 \rangle^{\frac{1}{2}-} \langle 
\sigma_2 \rangle^{\frac{1}{2}}} \, . $$
{\bf a.} $l+k \le 0$ . \\
Using (\ref{3.1}) we get
$$ S \le c\int \frac{|\widehat{v} \widehat{v_1} \widehat{v_2}| (\langle \sigma 
\rangle + \langle \sigma_1 \rangle + \langle \sigma_2 
\rangle)^{\frac{-k-l}{2}}}{|\xi_1|^{1-\epsilon} \langle 
\sigma \rangle^{\frac{1}{2}} \langle \sigma_1 \rangle^{\frac{1}{2}-} \langle 
\sigma_2 \rangle^{\frac{1}{2}}} \, . $$
Remark that $-k-l<1$ , so that Lemma \ref{LemmaB} can be applied with 
$m=1-\epsilon$ and gives $2(a+a_1+a_2)+m=k+l-4-\epsilon- \ge 3-\epsilon- > 2$, 
because $k\ge 0$ , $l\ge -1$ , thus the desired estimate follows. \\
{\bf b.} $l+k \ge 0$ . \\
$$ S \le c\int \frac{|\widehat{v} \widehat{v_1} 
\widehat{v_2}|}{|\xi_1|^{1-\epsilon} \langle 
\sigma \rangle^{\frac{1}{2}} \langle \sigma_1 \rangle^{\frac{1}{2}-} \langle 
\sigma_2 \rangle^{\frac{1}{2}}} \, . $$
Using Lemma \ref{LemmaB} with $m=1-\epsilon$ gives the desired result. \\
{\bf Case 3:} $|\xi_1| \le 1 $ , $|\xi_2| \le 1 $ and w.l.o.g. $\delta \le 
\epsilon$ . \\
$$ S \le c\int \frac{|\widehat{v} \widehat{v_1} \widehat{v_2}| 
|\xi_1|^{-1+\epsilon}}{|\xi_2|^{\epsilon - \delta} \langle 
\sigma \rangle^{\frac{1}{2}} \langle \sigma_1 \rangle^{\frac{1}{2}-} \langle 
\sigma_2 \rangle^{\frac{1}{2}}} \le c\int \frac{|\widehat{v} \widehat{v_1} 
\widehat{v_2}|}{|\xi_1|^{1 - \delta} \langle 
\sigma \rangle^{\frac{1}{2}} \langle \sigma_1 \rangle^{\frac{1}{2}-} \langle 
\sigma_2 \rangle^{\frac{1}{2}}}\, . $$
Using Lemma \ref{LemmaB} with $m=1-\delta$ gives the result. \\
{\bf Region B:} $\frac{1}{2}|\xi_2| \le |\xi_1| \le 2|\xi_2|$ ($\Rightarrow 
|\xi| \le 3|\xi_1|,3|\xi_2|$). \\
We have
$$ S \le c\int \frac{|\widehat{v} \widehat{v_1} \widehat{v_2}| \langle \xi 
\rangle^{-l-\delta} |\xi|^{\delta}}{|\xi_1| \langle 
\sigma \rangle^{\frac{1}{2}} \langle \sigma_1 \rangle^{\frac{1}{2}-} \langle 
\sigma_2 \rangle^{\frac{1}{2}}} \, . $$
{\bf Case 1:} $|\xi| \le 1$ . \\
$$ S \le c\int \frac{|\widehat{v} \widehat{v_1} 
\widehat{v_2}|}{|\xi_1|^{1-\delta} \langle 
\sigma \rangle^{\frac{1}{2}} \langle \sigma_1 \rangle^{\frac{1}{2}-} \langle 
\sigma_2 \rangle^{\frac{1}{2}}} \, . $$
This can easily be handled by Lemma \ref{LemmaB} with $m=1-\delta$ . \\
{\bf Case 2:} $|\xi| \ge 1 $ ($\Rightarrow |\xi_1| \ge \frac{1}{3}$).
$$ S \le c\int \frac{|\widehat{v} \widehat{v_1} \widehat{v_2}|\langle \xi 
\rangle^{-l}}{\langle \xi_1 \rangle \langle 
\sigma \rangle^{\frac{1}{2}} \langle \sigma_1 \rangle^{\frac{1}{2}-} \langle 
\sigma_2 \rangle^{\frac{1}{2}}} \, . $$
{\bf a.} $l \le 0$ . \\
$$ S \le c\int \frac{|\widehat{v} \widehat{v_1} \widehat{v_2}|}{\langle \xi_1 
\rangle^{1+l} \langle 
\sigma \rangle^{\frac{1}{2}} \langle \sigma_1 \rangle^{\frac{1}{2}-} \langle 
\sigma_2 \rangle^{\frac{1}{2}}} \, . $$ 
Because $1+l \ge 0$ this can easily be handled by Lemma \ref{LemmaA} or Lemma 
\ref{LemmaB}. \\
{\bf b.} $l \ge 0$ . \\
We get 
$$ S \le c\int \frac{|\widehat{v} \widehat{v_1} \widehat{v_2}|}{\langle \xi_1 
\rangle \langle 
\sigma \rangle^{\frac{1}{2}} \langle \sigma_1 \rangle^{\frac{1}{2}-} \langle 
\sigma_2 \rangle^{\frac{1}{2}}} \, , $$ 
which can be treated by Lemma \ref{LemmaA}. \\ 
{\bf Region C:} $|\xi_1| \ge 2|\xi_2|$ ( $\Rightarrow |\xi| \sim |\xi_1|$). \\
We get 
$$ S \le c\int \frac{|\widehat{v} \widehat{v_1} \widehat{v_2}| 
|\xi_1|^{-1+\epsilon +\delta} \langle \xi_1 
\rangle^{k-\epsilon-l-\delta}}{|\xi_2|^{\epsilon} \langle \xi_2 
\rangle^{k-\epsilon} \langle 
\sigma \rangle^{\frac{1}{2}} \langle \sigma_1 \rangle^{\frac{1}{2}-} \langle 
\sigma_2 \rangle^{\frac{1}{2}}} \, . $$ 
{\bf Case 1:} $|\xi_1| \ge 1 $ , $ |\xi_2| \ge 1 $ . \\
This case can be handled like the 3-dimensional case in Lemma \ref{LemmaC}. \\
{\bf Case 2:} $|\xi_1| \ge 1 $ , $ |\xi_2| \le 1 $ . \\ 
We have by (\ref{3.1}) : 
$$ S \le c\int \frac{|\widehat{v} \widehat{v_1} \widehat{v_2}| \langle \xi_1 
\rangle^{k-l-1}}{|\xi_2|^{\epsilon} \langle 
\sigma \rangle^{\frac{1}{2}} \langle \sigma_1 \rangle^{\frac{1}{2}-} \langle 
\sigma_2 \rangle^{\frac{1}{2}}} \le c\int \frac{|\widehat{v} \widehat{v_1} 
\widehat{v_2}| (\langle \sigma \rangle + \langle \sigma_1 \rangle + \langle 
\sigma_2 \rangle)^{\frac{k-l-1}{2}}}{|\xi_2|^{\epsilon} \langle 
\sigma \rangle^{\frac{1}{2}} \langle \sigma_1 \rangle^{\frac{1}{2}-} \langle 
\sigma_2 \rangle^{\frac{1}{2}}} \, .$$  
Because $k>l+2$ we can apply Lemma \ref{LemmaB} with $m=\epsilon$ and compute 
$2(a+a_1+a_2)+m=2(\frac{1}{2}+\frac{1}{2}+\frac{1}{2}-\frac{k-l-1}{2})+\epsilon- 
> 2+\epsilon-$ , so that the claimed estimate follows. \\
{\bf Case 3:} $|\xi_1| \le 1$ , $ |\xi_2| \le 1$ and w.l.o.g. $\delta \le 
1-\epsilon$ . \\
$$ S \le c\int \frac{|\widehat{v} \widehat{v_1} \widehat{v_2}|  
|\xi_1|^{-1+\epsilon+\delta}}{|\xi_2|^{\epsilon} \langle 
\sigma \rangle^{\frac{1}{2}} \langle \sigma_1 \rangle^{\frac{1}{2}-} \langle 
\sigma_2 \rangle^{\frac{1}{2}}} \le c\int \frac{|\widehat{v} \widehat{v_1} 
\widehat{v_2}|}{|\xi_2|^{1-\delta} \langle 
\sigma \rangle^{\frac{1}{2}} \langle \sigma_1 \rangle^{\frac{1}{2}-} \langle 
\sigma_2 \rangle^{\frac{1}{2}}}\, . $$
An application of Lemma \ref{LemmaB} with $m=1-\delta$ gives the desired 
estimate.
\begin{lemma}
\label{LemmaF}
Let $n=2$. Under the assumptions of Lemma \ref{LemmaD} we have
$$\|B^{2-\delta} (D\bar{\varphi} D\varphi)\|_{X_{\pm}^{l+\delta,-\frac{1}{2}+}} 
\le c T^{\Theta} \|B^{\epsilon} D\varphi\|_{X^{k-\epsilon,\frac{1}{2}}}^2 $$
with $\Theta > 0$ for $0< \delta < 1$ , $0 < \epsilon < 1$ .
\end{lemma} 
{\bf Proof:} Arguing as in Lemma \ref{LemmaD} we have to show 
$$
W:=\left| \int \frac{\widehat{v} \widehat{v_1} \widehat{v_2} \langle \xi 
\rangle^{l+\delta} 
 |\xi|^{2-\delta}}{|\xi_1|^{\epsilon}\langle \xi_1 \rangle^{k-\epsilon} 
|\xi_2|^{\epsilon} \langle \xi_2 \rangle^{k-\epsilon} \langle 
\sigma \rangle^{\frac{1}{2}-} \langle \sigma_1 \rangle^{\frac{1}{2}} \langle 
\sigma_2 \rangle^{\frac{1}{2}}} \right| \le cT^{\Theta} \|v\|_{L^2_{xt}} 
\|v_1\|_{L^2_{xt}} \|v_2\|_{L^2_{xt}} \, .
$$
{\bf Region A:} $\frac{|\xi_2|}{2} \le |\xi_1| \le 2|\xi_2|$ ($\Rightarrow |\xi| 
\le 3|\xi_1|,3|\xi_2|$). \\
{\bf Case 1:} $|\xi_1| \ge 1$ ($\Rightarrow |\xi_2| \ge \frac{1}{2}$) . \\
Using the assumption $k \ge \frac{l+2}{2}$ we get 
\begin{eqnarray*}
W \le c\int \frac{|\widehat{v} \widehat{v_1} \widehat{v_2}| \langle \xi 
\rangle^{l+2}}{\langle \xi_1 \rangle^k \langle \xi_2 \rangle^k \langle 
\sigma \rangle^{\frac{1}{2}-} \langle \sigma_1 \rangle^{\frac{1}{2}} \langle 
\sigma_2 \rangle^{\frac{1}{2}}} & \le & c\int \frac{|\widehat{v} \widehat{v_1} 
\widehat{v_2}| \langle \xi \rangle^{l+2-2k}}{\langle 
\sigma \rangle^{\frac{1}{2}-} \langle \sigma_1 \rangle^{\frac{1}{2}} \langle 
\sigma_2 \rangle^{\frac{1}{2}}} \\
& \le & c\int \frac{|\widehat{v} \widehat{v_1} \widehat{v_2}|}{\langle 
\sigma \rangle^{\frac{1}{2}-} \langle \sigma_1 \rangle^{\frac{1}{2}} \langle 
\sigma_2 \rangle^{\frac{1}{2}}} \, .
\end{eqnarray*}
Lemma \ref{LemmaB} gives the claimed estimate. \\
{\bf Case 2:} $|\xi_1| \le 1$ ($\Rightarrow |\xi_2| \le 2 \Rightarrow |\xi| \le 
3$) . \\
Using $2-\delta-\epsilon >0$ and $|\xi| \le 3$ we get the estimate 
\begin{eqnarray*}
W \le c\int \frac{|\widehat{v} \widehat{v_1} \widehat{v_2}| 
|\xi|^{2-\delta}}{|\xi_1|^{\epsilon} |\xi_2|^{\epsilon} \langle 
\sigma \rangle^{\frac{1}{2}-} \langle \sigma_1 \rangle^{\frac{1}{2}} \langle 
\sigma_2 \rangle^{\frac{1}{2}}} & \le & c\int \frac{|\widehat{v} \widehat{v_1} 
\widehat{v_2}| |\xi|^{2-\delta-\epsilon}}{ 
|\xi_2|^{\epsilon} \langle \sigma \rangle^{\frac{1}{2}-} \langle \sigma_1 
\rangle^{\frac{1}{2}} \langle \sigma_2 \rangle^{\frac{1}{2}}} \\
& \le & c\int \frac{|\widehat{v} \widehat{v_1} 
\widehat{v_2}|}{|\xi_2|^{\epsilon} \langle 
\sigma \rangle^{\frac{1}{2}-} \langle \sigma_1 \rangle^{\frac{1}{2}} \langle 
\sigma_2 \rangle^{\frac{1}{2}}} \, .
\end{eqnarray*} 
Lemma \ref{LemmaB} gives the claimed estimate. \\
{\bf Region B:} $|\xi_1| \ge 2|\xi_2|$ ($ \Rightarrow |\xi| \sim |\xi_1|$) (and 
similarly $|\xi_2| \ge 2|\xi_1|$). \\
We get 
$$
W \hspace{-0.1em}\le \hspace{-0.1em}c \int \frac{|\widehat{v} \widehat{v_1} 
\widehat{v_2}| \langle \xi_1 \rangle^{l+\delta} 
|\xi_1|^{2-\delta-\epsilon}}{\langle \xi_1 \rangle^{k-\epsilon} 
|\xi_2|^{\epsilon} \langle \xi_2 \rangle^{k-\epsilon} \langle 
\sigma \rangle^{\frac{1}{2}-} \langle \sigma_1 \rangle^{\frac{1}{2}} \langle 
\sigma_2 \rangle^{\frac{1}{2}}}
\le c \int \frac{|\widehat{v} \widehat{v_1} \widehat{v_2}| \langle \xi_1 
\rangle^{l+2-k}}{|\xi_2|^{\epsilon} \langle \xi_2 \rangle^{k-\epsilon} \langle 
\sigma \rangle^{\frac{1}{2}-} \langle \sigma_1 \rangle^{\frac{1}{2}} \langle 
\sigma_2 \rangle^{\frac{1}{2}}}. 
$$ 
{\bf Case 1:} $|\xi_2| \ge 1$ . \\
$$
W \le c \int \frac{|\widehat{v} \widehat{v_1} \widehat{v_2}| \langle \xi_1 
\rangle^{l+2-k}}{\langle \xi_2 \rangle^k \langle 
\sigma \rangle^{\frac{1}{2}-} \langle \sigma_1 \rangle^{\frac{1}{2}} \langle 
\sigma_2 \rangle^{\frac{1}{2}}} \, . $$
This is exactly the integral treated in Lemma \ref{LemmaD} in the case $n=3$. \\ 
{\bf Case 2:} $|\xi_2| \le 1$ . \\
Assuming w.l.o.g. $k \le l+2$ and using (\ref{3.1}) we get the estimate 
$$
W \le c \int \frac{|\widehat{v} \widehat{v_1} \widehat{v_2}| \langle \xi_1 
\rangle^{l+2-k}}{|\xi_2|^{\epsilon} \langle 
\sigma \rangle^{\frac{1}{2}-} \langle \sigma_1 \rangle^{\frac{1}{2}} \langle 
\sigma_2 \rangle^{\frac{1}{2}}} 
\le  \int \frac{|\widehat{v} \widehat{v_1} \widehat{v_2}| (\langle \sigma 
\rangle + \langle \sigma_1 \rangle + \langle \sigma_2 
\rangle)^{\frac{l+2-k}{2}}}{|\xi_2|^{\epsilon} \langle 
\sigma \rangle^{\frac{1}{2}-} \langle \sigma_1 \rangle^{\frac{1}{2}} \langle 
\sigma_2 \rangle^{\frac{1}{2}}} \, . $$ 
The exponents in the denominator are nonnegative, because $k>l+1$ . Thus we 
apply Lemma \ref{LemmaB} with e.g. $a=\frac{1}{2} - \frac{l+2-k}{2}- > 0$ , 
$a_1=a_2=\frac{1}{2}$ , $m=\epsilon$ , so that $2(a+a_1+a_2)+m > 2+\epsilon > 2$ 
.  

The following variant of Lemma \ref{LemmaF} is also true:
\begin{lemma}
\label{Lemma tildeF}
Let $n=2$ . Under the assumptions of Lemma \ref{LemmaD} we have
$$ \|D\bar{\varphi} D\varphi\|_{\dot{X}_{\pm}^{l+2,-\frac{1}{2}+}} \sim  
\|B^{2-\delta}D\bar{\varphi} D\varphi\|_{\dot{X}_{\pm}^{l+\delta,-\frac{1}{2}+}} 
\le cT^{\Theta} \|B^{\epsilon} D\varphi\|_{X^{k-\epsilon,\frac{1}{2}}}^2 $$
with $\Theta >0$ for $0<\epsilon <1$ .
\end{lemma}  
{\bf Proof:} The proof of Lemma \ref{LemmaF} is modified as follows. We have to 
estimate 
$$
W:=\left| \int \frac{\widehat{v} \widehat{v_1} \widehat{v_2} 
|\xi|^{l+2}}{|\xi_1|^{\epsilon}\langle \xi_1 \rangle^{k-\epsilon} 
|\xi_2|^{\epsilon} \langle \xi_2 \rangle^{k-\epsilon} \langle 
\sigma \rangle^{\frac{1}{2}-} \langle \sigma_1 \rangle^{\frac{1}{2}} \langle 
\sigma_2 \rangle^{\frac{1}{2}}} \right| \, . $$
{\bf Region A:} $\frac{|\xi_2|}{2} \le |\xi_1| \le 2|\xi_2|$ ($\Rightarrow |\xi| 
\le 3|\xi_1|,3|\xi_2|$). \\
{\bf Case 1:} $|\xi_1| \ge 1$ ($\Rightarrow |\xi_2| \ge \frac{1}{2}$) . \\  
This case is treated exactly as in Lemma \ref{LemmaF}. \\
{\bf Case 2:} $|\xi_1| \le 1$ ($\Rightarrow |\xi_2| \le 2 \Rightarrow |\xi| \le 
3$) . \\
Using $l+2-\epsilon >0$ and $|\xi| \le 3$ we get the bound 
\begin{eqnarray*}
W \le c\int \frac{|\widehat{v} \widehat{v_1} \widehat{v_2}| \, 
|\xi|^{l+2}}{|\xi_1|^{\epsilon} |\xi_2|^{\epsilon} \langle 
\sigma \rangle^{\frac{1}{2}-} \langle \sigma_1 \rangle^{\frac{1}{2}} \langle 
\sigma_2 \rangle^{\frac{1}{2}}}
& \le & c\int \frac{|\widehat{v} \widehat{v_1} \widehat{v_2}| \, 
|\xi|^{l+2-\epsilon}}{|\xi_2|^{\epsilon} \langle 
\sigma \rangle^{\frac{1}{2}-} \langle \sigma_1 \rangle^{\frac{1}{2}} \langle 
\sigma_2 \rangle^{\frac{1}{2}}} \\
& \le & c\int \frac{|\widehat{v} \widehat{v_1} 
\widehat{v_2}|}{|\xi_2|^{\epsilon} \langle 
\sigma \rangle^{\frac{1}{2}-} \langle \sigma_1 \rangle^{\frac{1}{2}} \langle 
\sigma_2 \rangle^{\frac{1}{2}}} \, , 
\end{eqnarray*} 
which can be estimated by Lemma \ref{LemmaB}. \\
{\bf Region B:} $|\xi_1| \ge 2|\xi_2|$ ($ \Rightarrow |\xi| \sim |\xi_1|$) (and 
similarly $|\xi_2| \ge 2|\xi_1|$). \\
Using $l+2-\epsilon >0$ we get the bound 
$$
W \hspace{-0.1em}\le \hspace{-0.1em}c \int \frac{|\widehat{v} \widehat{v_1} 
\widehat{v_2}| |\xi_1|^{l+2-\epsilon}}{\langle \xi_1 \rangle^{k-\epsilon} 
|\xi_2|^{\epsilon} \langle \xi_2 \rangle^{k-\epsilon} \langle 
\sigma \rangle^{\frac{1}{2}-} \langle \sigma_1 \rangle^{\frac{1}{2}} \langle 
\sigma_2 \rangle^{\frac{1}{2}}}
\le c \int \frac{|\widehat{v} \widehat{v_1} \widehat{v_2}| \langle \xi_1 
\rangle^{l+2-k}}{|\xi_2|^{\epsilon} \langle \xi_2 \rangle^{k-\epsilon} \langle 
\sigma \rangle^{\frac{1}{2}-} \langle \sigma_1 \rangle^{\frac{1}{2}} \langle 
\sigma_2 \rangle^{\frac{1}{2}}} . $$ 
This is exactly the integral treated in the proof of Lemma \ref{LemmaF}, Region 
B. Thus the claimed estimate follows.

In order to treat the limiting cases $k=l+1$ and $k=l+2$ we also need the 
following results:
\begin{lemma}
\label{LemmaE'}
Let $n=2$ , $l \ge -1$ , $l+1 \le k \le l+2$ , and let $\varphi,\chi$ be 
supported in $\{|t| \le cT\}$ . Then the following estimate holds:
$$ \|B^{-1+\epsilon} (D\varphi D\chi)\|_{X^{k-\epsilon,-\frac{1}{2}}} \le c 
T^{\Theta} \|B^{\epsilon} D\varphi\|_{X^{k-\epsilon,\frac{1}{2}}} \|B^{-\delta} 
D\chi\|_{X_{\pm}^{l+\delta,\frac{1}{2}}} $$
with $\Theta >0$ for $0 < \epsilon < 1 $ , $ \delta > 0 $ .
\end{lemma}
{\bf Remark:} For $l<0$ we can replace $\|B^{-\delta} 
D\chi\|_{X_{\pm}^{l+\delta,\frac{1}{2}}} $ by 
$\|\chi\|_{\dot{X}_{\pm}^{l+1,\frac{1}{2}}}$ . \\
{\bf Proof:} We repeat the proof of Lemma \ref{LemmaE} replacing $\langle 
\sigma_1 \rangle^{\frac{1}{2}-}$ by $\langle \sigma_1 \rangle^{\frac{1}{2}}$ . 
We only have to remark that the limit case $k=l+2$ is allowed in Region C, Case 
1 and Case 2, because the power of the $\sigma$ - modules in the denominator 
remains nonnegative in this case.
\begin{lemma}
\label{LemmaF'}
Let $n=2$ , $l \ge -1$ , $k \ge \frac{l+2}{2}$ , $k = l+1 $ and supp $\varphi 
\subset \{|t| \le cT\}$ . Then
$$ \|B^{2-\delta} (D\bar{\varphi} D\varphi)\|_{X_{\pm}^{l+\delta,-\frac{1}{2}}} 
\le c T^{\Theta} \|B^{\epsilon} D\varphi\|_{X^{k-\epsilon,\frac{1}{2}}}^2 $$
with $\Theta >0$ for $0<\delta <1$ , $0< \epsilon < 1$ .
\end{lemma}
{\bf Proof:} We repeat the proof of Lemma \ref{LemmaF} with $\langle \sigma 
\rangle^{\frac{1}{2}-}$ replaced by  $\langle \sigma \rangle^{\frac{1}{2}}$ . 
The condition $k<l+1$ was only used in Region B, Cases 1 and 2 to produce 
nonnegative exponents of the $\sigma$ - modules in the denominator, which is 
satisfied now also for $k=l+1$ . \\
{\bf Remark:} The estimate of Lemma \ref{Lemma tildeF} remains true for $k=l+1$ 
in the following form:
$$ \|D\bar{\varphi} D\varphi\|_{\dot{X}_{\pm}^{l+2,-\frac{1}{2}}} \le 
cT^{\Theta} \|B^{\epsilon} D\varphi\|_{X^{k-\epsilon,\frac{1}{2}}}^2 $$
with $\Theta >0$ for $0<\epsilon <1$ .\\ This follows similarly as Lemma 
\ref{LemmaF'}. 
\begin{lemma}
\label{LemmaE''}
Assume $n=2$ , $l \ge -1$ , $k=l+2$, and let $\varphi,\chi$ be supported in 
$\{|t|\le cT\}$. Then
$$ \|B^{-1+\epsilon} (D\varphi D\chi)\|_{Y^{k-\epsilon}} \le c T^{\Theta} 
\|B^{\epsilon} D\varphi\|_{X^{k-\epsilon,\frac{1}{2}}} \|B^{-\delta} 
D\chi\|_{X_{\pm}^{l+\delta,\frac{1}{2}}} $$
with $\Theta >0$ for $0 < \epsilon < 1 $ , $ \delta > 0 $ .
\end{lemma} 
{\bf Remark:} For $l<0$ we can replace $\|B^{-\delta} 
D\chi\|_{X_{\pm}^{l+\delta,\frac{1}{2}}}$ by 
$\|\chi\|_{\dot{X}^{l+1,\frac{1}{2}}}$ . \\
{\bf Proof:} Arguing as in the proof of Lemma \ref{LemmaE} we now have to give 
the following estimate (cf. the proof of Lemma \ref{LemmaC''}):
$$
\tilde{S}:= \int \frac{|\widehat{v} \widehat{w_1} \widehat{v_2}| \, 
|\xi_1|^{-1+\epsilon} 
\langle \xi_1 \rangle^{k-\epsilon}}{|\xi|^{-\delta}\langle \xi 
\rangle^{l+\delta} |\xi_2|^{\epsilon} \langle \xi_2 \rangle^{k-\epsilon} \langle 
\sigma \rangle^{\frac{1}{2}} \langle \sigma_1 \rangle \langle 
\sigma_2 \rangle^{\frac{1}{2}}} \le cT^{\Theta} \|v\|_{L^2_{xt}} 
\|w_1\|_{L^2_x} \|v_2\|_{L^2_{xt}} \, .
$$
The only case where the strict inequality $k<l+2$ was used in the proof of Lemma 
\ref{LemmaE} was Region C, Case 1 and 2. In all other regions we define 
$\widehat{v_1} :=\langle \sigma_1 \rangle^{-\frac{1}{2}-} \widehat{w_1},$ so 
that $\|v_1\|_{L^2_{xt}} \le c \|w_1\|_{L_x^2}$ , and $\tilde{S}$ reads  as 
follows: 
$$
\tilde{S} = \int \frac{|\widehat{v} \widehat{v_1} \widehat{v_2}| \, 
|\xi_1|^{-1+\epsilon} 
\langle \xi_1 \rangle^{k-\epsilon}}{|\xi|^{-\delta}\langle \xi 
\rangle^{l+\delta} |\xi_2|^{\epsilon} \langle \xi_2 \rangle^{k-\epsilon} \langle 
\sigma \rangle^{\frac{1}{2}} \langle \sigma_1 \rangle^{\frac{1}{2}-} \langle 
\sigma_2 \rangle^{\frac{1}{2}}} \, . $$
This is exactly the integral treated in the proof of Lemma \ref{LemmaE}, so that 
the result in these regions follows. It remains to consider Region C, Case 1 and 
2 in the proof of Lemma \ref{LemmaE}. Similarly as there we get in Region C, 
Case 1 (with $k=l+2$): 
$$
\tilde{S} \le c \int \frac{|\widehat{v} \widehat{w_1} \widehat{v_2}|  
\langle \xi_1 \rangle}{\langle \xi_2 \rangle^k \langle 
\sigma \rangle^{\frac{1}{2}} \langle \sigma_1 \rangle \langle 
\sigma_2 \rangle^{\frac{1}{2}}} \, . $$ 
This integral was already treated in the proof of Lemma \ref{LemmaC''}. In 
Region C, Case 2 by use of (\ref{3.2}) we arrive at
$$
\tilde{S} \le c \int \frac{|\widehat{v} \widehat{w_1} \widehat{v_2}|  
\langle \xi_1 \rangle}{|\xi_2|^{\epsilon} \langle 
\sigma \rangle^{\frac{1}{2}} \langle \sigma_1 \rangle \langle 
\sigma_2 \rangle^{\frac{1}{2}}} \le c \int \frac{|\widehat{v} \widehat{w_1} 
\widehat{v_2}| (\langle \sigma \rangle + \langle \sigma_2 \rangle + \langle 
\sigma_1 \rangle \phi_{\{c_1|\sigma_1| \le |\xi_1|^2 \le 
c_2|\sigma_1|\}})^{\frac{1}{2}}}{|\xi_2|^{\epsilon} \langle 
\sigma \rangle^{\frac{1}{2}} \langle \sigma_1 \rangle \langle 
\sigma_2 \rangle^{\frac{1}{2}}} \, .$$ 
The two terms coming from $\langle \sigma \rangle$ and $\langle \sigma_2 
\rangle$ in the numerator are treated by defining $\widehat{v_1}$ as before by 
Lemma \ref{LemmaB} with e.g. $a=0$ , $a_1= \frac{1}{2}-$ , $a_2=\frac{1}{2}$ , 
$m=\epsilon$ , so that $2(a+a_1+a_2)+m=2+\epsilon- >2$ , whereas the term coming 
from $\langle \sigma_1 \rangle$ is treated by defining $\widehat{v_1}:=\langle 
\sigma_1 \rangle^{-\frac{1}{2}} \widehat{w_1} \phi_{\{c_1 |\sigma_1| \le 
|\xi_1|^2 \le c_2|\sigma_1|\}}$ . so that $\|v_1\|_{L^2_{xt}} \le c 
\|w_1\|_{L_x^2}$ 
. Thus we are left with 
$$\int \frac{|\widehat{v} \widehat{v_1} \widehat{v_2}|}{|\xi_2|^{\epsilon} 
\langle 
\sigma \rangle^{\frac{1}{2}} \langle \sigma_2 \rangle^{\frac{1}{2}}} \, , $$
which can be handled by Lemma \ref{LemmaB}. 

Finally we get
\begin{lemma}
\label{LemmaF''}
Let $n=2$ , $l \ge -1$ , $k \ge \frac{l+2}{2}$ , $k = l+1 $ and supp $\varphi 
\subset \{|t| \le cT\}$ . Then
$$ \|B^{2-\delta} (D\bar{\varphi} D\varphi)\|_{Y_{\pm}^{l+\delta}} \le c 
T^{\Theta} \|B^{\epsilon} D\varphi\|_{X^{k-\epsilon,\frac{1}{2}}}^2 $$
with $\Theta >0$ for $0<\delta <1$ , $0< \epsilon < 1$ .
\end{lemma} 
{\bf Proof:} We follow the proof of Lemma \ref{LemmaF} and have to show 
$$
\tilde{W}:= \int \frac{|\widehat{w} \widehat{v_1} \widehat{v_2}| \langle \xi 
\rangle^{l+\delta} 
 |\xi|^{2-\delta}}{|\xi_1|^{\epsilon}\langle \xi_1 \rangle^{k-\epsilon} 
|\xi_2|^{\epsilon} \langle \xi_2 \rangle^{k-\epsilon} \langle 
\sigma \rangle \langle \sigma_1 \rangle^{\frac{1}{2}} \langle 
\sigma_2 \rangle^{\frac{1}{2}}} \le cT^{\Theta} \|w\|_{L^2_x} 
\|v_1\|_{L^2_{xt}} \|v_2\|_{L^2_{xt}} \, .
$$ 
In Region A, Case 1 of the proof of Lemma \ref{LemmaF} we define 
$\widehat{v}:=\langle \sigma \rangle^{-\frac{1}{2}-} \widehat{w}$ , so that 
$\|v\|_{L^2_{xt}} \le c \|w\|_{L_x^2}$ , and we get as in Lemma \ref{LemmaF} the 
estimate 
$$
\tilde{W} \le c \int \frac{|\widehat{w} \widehat{v_1} \widehat{v_2}|}{\langle 
\sigma \rangle \langle \sigma_1 \rangle^{\frac{1}{2}} \langle 
\sigma_2 \rangle^{\frac{1}{2}}} \le c \int \frac{|\widehat{v} \widehat{v_1} 
\widehat{v_2}|}{\langle 
\sigma \rangle^{\frac{1}{2}-} \langle \sigma_1 \rangle^{\frac{1}{2}} \langle 
\sigma_2 \rangle^{\frac{1}{2}}} \, , $$
which can easily be handled by Lemma \ref{LemmaB}. \\
Similarly, in Region A, Case 2 we arrive at 
$$ \tilde{W} \le c \int \frac{|\widehat{v} \widehat{v_1} 
\widehat{v_2}|}{|\xi_1|^{\epsilon} \langle 
\sigma \rangle^{\frac{1}{2}-} \langle \sigma_1 \rangle^{\frac{1}{2}} \langle 
\sigma_2 \rangle^{\frac{1}{2}}} \, , $$ 
which can be controlled by Lemma \ref{LemmaB} again. \\
In Region B, Case 1 we get for $k=l+1$ using (\ref{3.3}): 
$$ \tilde{W} \le c \int \frac{|\widehat{w} \widehat{v_1} \widehat{v_2}| \langle 
\xi_1 \rangle}{\langle \xi_2 \rangle^k \langle 
\sigma \rangle \langle \sigma_1 \rangle^{\frac{1}{2}} \langle 
\sigma_2 \rangle^{\frac{1}{2}}} \le c \int \frac{|\widehat{w} \widehat{v_1} 
\widehat{v_2}| (\langle \sigma_1 \rangle + \langle \sigma_2 \rangle + \langle 
\sigma \rangle \phi_{\{c_1|\sigma| \le |\xi|^2 \le 
c_2|\sigma|\}})^{\frac{1}{2}}}{\langle \xi_2 \rangle^k \langle 
\sigma \rangle \langle \sigma_1 \rangle^{\frac{1}{2}} \langle 
\sigma_2 \rangle^{\frac{1}{2}}} \, . $$
The two terms coming from $\langle \sigma_1 \rangle$ and $\langle \sigma_2 
\rangle$ in the numerator are treated by defining $\widehat{v}$ as before and 
using Lemma \ref{LemmaA}, whereas the term coming from $\langle \sigma \rangle$ 
is treated by defining $\widehat{v}:= \langle \sigma \rangle^{-\frac{1}{2}} 
\widehat{w} \phi_{\{c_1|\sigma| \le |\xi|^2 \le c_2|\sigma| \}}$ 
, so that $\|v\|_{L^2_{xt}} \le c \|w\|_{L^2_x} $ , leading to 
$$ \int \frac{|\widehat{v} \widehat{v_1} \widehat{v_2}|}{\langle \xi_2 \rangle^k 
 \langle \sigma_1 \rangle^{\frac{1}{2}} \langle 
\sigma_2 \rangle^{\frac{1}{2}}} \, , $$
which again can be handled by Lemma \ref{LemmaA} (remark that $k \ge 
\frac{1}{2}$). \\
In Region B, Case 2 we arrive at the corresponding integrals where $\langle 
\xi_2 \rangle^k$ is replaced by $|\xi_2|^{\epsilon}$ . This can be treated by 
use of Lemma \ref{LemmaB}. \\
{\bf Remark:} The following variant of Lemma \ref{LemmaF''} is also true, as 
follows similarly from the proof of Lemma \ref{Lemma tildeF}: \\ 
Let $n=2$ , $l \ge -1$ , $k \ge \frac{l+2}{2}$ , $k = l+1 $ and supp $\varphi 
\subset \{|t| \le cT\}$ . Then
$$ \| (D\bar{\varphi} D\varphi)\|_{\dot{Y}_{\pm}^{l+2}} \le c T^{\Theta} 
\|B^{\epsilon} D\varphi\|_{X^{k-\epsilon,\frac{1}{2}}}^2 $$
with $\Theta >0$ for $0< \epsilon < 1$ .

These results can now be used to prove a local existence and uniqueness result 
as in the 3+1-dimensional case.
\begin{theorem}
\label{Theorem5}
In space dimension $n=2$ assume $l \ge -1$ , $l+1 \le k \le l+2, $  $ k \ge 
\frac{l+2}{2}$ , $0<\epsilon,\delta <1$, and 
$$B^{1+\epsilon}\varphi_0 \in H^{k-\epsilon}({\bf R^2}) \, , \, 
B^{1-\delta}\chi_0 \in H^{l+\delta}({\bf R}^2) \, ,  \, B^{-\delta}\chi_1 \in 
H^{l+\delta}({\bf R}^2). $$ Then there exists $1 \ge   
T=T(\|B^{1+\epsilon}\varphi_0\|_{H^{k-\epsilon}},\|B^{1-\delta}\chi_0\|_{H^{l+
\delta}},\|B^{-\delta}\chi_1\|_{H^{l+\delta}}) > 0 $ , such that the problem 
(\ref{1.4}),(\ref{1.5}),(\ref{1.3}) has a unique solution $(\varphi,\chi)$ with 
$$B^{1+\epsilon}\varphi \in X^{k-\epsilon,b}[0,T]\quad , \quad B^{1-\delta}\chi 
\, , \, B^{-\delta}\chi_t \in X_+^{l+\delta,b_1}[0,T] + X_-^{l+\delta,b_1}[0,T] 
\, .$$ Here $b=\frac{1}{2}+$ , $b_1=\frac{1}{2}+$ , if $l+1<k<l+2$ , 
$b=\frac{1}{2}$ , $b_1=\frac{1}{2}+$ , if $ k=l+2$ , and $b=\frac{1}{2}+,$ 
$b_1=\frac{1}{2}$ , if $k=l+1$ . This solution satisfies $$B^{1+\epsilon}\varphi 
\in C^0([0,T],H^{k-\epsilon}({\bf R}^2)) \, , \, B^{1-\delta}\chi \, , \, 
B^{-\delta}\chi_t \in C^0([0,T],H^{l+\delta}({\bf R}^2)) . $$
If $l < 0$ we can replace $B^{1-\delta}\chi_0 \, , \, B^{-\delta}\chi_1 \in 
H^{l+\delta}$ by $\chi_0 \in \dot{H}^{l+1} $ , $\chi_1 \in \dot{H}^l$ , and 
$B^{1-\delta} \chi$ , $ B^{-\delta} \chi_t \in X_+^{l+\delta,b_1}[0,T] + 
X_-^{l+\delta,b_1}[0,T] $ by $\chi \in \dot{X}_+^{l+1,b_1}[0,T] + 
\dot{X}_-^{l+1,b_1}[0,T],$ $ \chi_t \in \dot{X}_+^{l,b_1}[0,T] + 
\dot{X}_-^{l,b_1}[0,T]$ , and we have $\chi \in C^0([0,T],\dot{H}^{l+1}({\bf 
R}^2)$ , $\chi_t \in C^0([0,T],\dot{H}^{l}({\bf R}^2)$ .
\end{theorem}  
{\bf Remark:} If this theorem would be true for $\epsilon =0$ , we would have 
local existence und uniqueness for data $B\varphi_0 \in H^1({\bf R}^2)$ , 
$\chi_0 \in L^2({\bf R}^2) $ , $ B^{-1} \chi_1 \in L^2({\bf R}^2)$ . Using the 
a-priori bounds for  $\|B\varphi\|_{H^1} + \|\chi\|_{L^2} + 
\|B^{-1}\chi_t\|_{L^2}$ under a smallness assumption on $\|B\varphi_0\|_{L^2}$ 
(cf. chapter 1) , this would imply global existence in these spaces under this 
smallness assumption.

\end{document}